\documentclass[11pt,reqno]{amsart}

\usepackage{amsmath,amssymb,amsthm}
\usepackage{graphicx}
\usepackage{booktabs}
\usepackage{mathrsfs}
\usepackage[margin=1in]{geometry}
\usepackage{hyperref}
\hypersetup{colorlinks=true,linkcolor=blue,citecolor=blue,urlcolor=blue,
  pdftitle={Sampled resolvents of C0-semigroups: aliasing and sharp rates on the Crandall-Pazy scale},
  pdfauthor={N. Castillo}}

\theoremstyle{plain}
\newtheorem{theorem}{Theorem}[section]
\newtheorem{proposition}[theorem]{Proposition}
\newtheorem{lemma}[theorem]{Lemma}
\newtheorem{corollary}[theorem]{Corollary}
\theoremstyle{definition}
\newtheorem{definition}[theorem]{Definition}

\theoremstyle{remark}
\newtheorem{remark}[theorem]{Remark}

\newcommand{\R}{\mathbb{R}}
\newcommand{\Z}{\mathbb{Z}}
\newcommand{\C}{\mathbb{C}}
\newcommand{\N}{\mathbb{N}}
\newcommand{\Hil}{\mathscr{H}}

\newcommand{\ip}[2]{\langle #1,#2\rangle}

\newcommand{\Qt}[1]{Q_{#1}}
\newcommand{\Qtr}[1]{Q^{\natural}_{#1}}
\newcommand{\F}[2]{F_{#1}^{#2}}
\newcommand{\eps}{\varepsilon}
\DeclareMathOperator{\Lis}{Li}
\DeclareMathOperator{\supp}{supp}
\DeclareMathOperator{\diam}{diam}
\renewcommand{\Re}{\operatorname{Re}}
\renewcommand{\Im}{\operatorname{Im}}

\begin{document}

\title[Sampled resolvents and the Crandall--Pazy scale]
{Sampled resolvents of $C_0$-semigroups:\\
aliasing and sharp rates on the Crandall--Pazy scale}

\author{N. Castillo}
\address{Independent researcher, Columbus, OH, USA}
\email{nickandpreeti@gmail.com}

\date{\today}


\begin{abstract}
For a $C_0$-semigroup $(S_t)$ with generator $B$, sampling the Laplace integral
$(\lambda-B)^{-1}=\int_0^\infty e^{-\lambda t}S_t\,dt$ on the grid $T\N$ gives
$Q_T=T(I-e^{-\lambda T}S_T)^{-1}$, the impulse-invariant discretization of the
resolvent. The motivation is a dyadic series of Castillo, Costin and Costin for the
resolvent of a self-adjoint operator: its partial sums are exactly $Q_{2^{-N}}$, so
the series is a telescoped rectangle rule. By Poisson summation, on $D(B)$ the error of
$Q_T$ is the scalar $T/2$ plus the symmetric limit of the resolvents at the alias points
$\lambda-2\pi in/T$, $n\ne0$, and the dyadic levels remove these aliases according to
the $2$-adic valuation of $n$. Our main result is a sharp scale of rates for the
trapezoid correction $Q_T-\frac T2I$. If $(S_t)$ belongs to the Crandall--Pazy class
with parameter $\alpha\in(0,1]$, that is $\|BS_t\|=O(t^{-1/\alpha})$ as $t\to0$, the
error is $O(T^\alpha)$ in operator norm, $O(T^{\alpha+\theta})$ on the interpolation
space of order $\theta$ between $X$ and $D(B)$, and $O(T^{1+\alpha})$ on $D(B)$, the
last requiring only resolvent decay of order $\alpha$ along the line
$\Re z=\Re\lambda$. This interpolates between unitary groups with unbounded generator
($\alpha=0$, no convergence in norm) and analytic semigroups ($\alpha=1$), and every exponent is
attained. For smooth data an Euler--Maclaurin expansion justifies Romberg
extrapolation. For self-adjoint generators we give closed forms for the contribution
of each dyadic level to Stone's formula, and we translate the results into an
explicit error bound and step-size rule for Green functions computed from sampled
autocorrelations. Experiments, including a structurally damped string whose semigroup
is of Gevrey class, confirm the rates.
\end{abstract}

\maketitle

%% =====================================================================
\section{Introduction}\label{sec:intro}

Let $(S_t)_{t\ge0}$ be a $C_0$-semigroup with generator $B$ on a Banach space $X$, and
let $\Re\lambda$ exceed its growth bound. Sampling the Laplace integral
$(\lambda-B)^{-1}=\int_0^\infty e^{-\lambda t}S_t\,dt$ at the times $mT$ gives
\[
\Qt{T}:=T\sum_{m\ge0}e^{-\lambda mT}S_{mT}=T\,(I-e^{-\lambda T}S_T)^{-1},
\]
the left-endpoint rectangle rule with step $T$. In signal-processing terms $\Qt{T}$ is
the impulse-invariant discretization of the transfer function
$(\lambda-B)^{-1}$~\cite{OS}. This paper studies the error
$\Qt{T}-(\lambda-B)^{-1}$: its exact structure, and its size in terms of the
smoothing properties of the semigroup.

The motivation is the following identity. For a self-adjoint operator $A$ on a
Hilbert space with unitary group $U_t=e^{-itA}$, Castillo, Costin and
Costin~\cite[Proposition~22(i)]{CCC} proved that for $\lambda>0$
\begin{equation}\label{eq:CCC}
(A-i\lambda)^{-1}
= i\,(1-e^{-\lambda}U_1)^{-1}
   - i\sum_{k=1}^{\infty}2^{-k}\bigl(1+e^{-\lambda 2^{-k}}U_{2^{-k}}\bigr)^{-1},
\end{equation}
with convergence in the strong operator topology, and in~\cite{CastilloV1} it was used
to derive dyadic series for lattice Green functions, fundamental solutions, spectral
measures and spectral zeta functions. With $S_t=U_t$ and $B=-iA$, the $k$-th term
of~\eqref{eq:CCC} is exactly the difference of the rectangle rules with steps
$2^{1-k}$ and $2^{-k}$: splitting the nodes of the finer rule into even and odd
indices gives $\Qt{2T}=\Qt{T}+T(1+e^{-\lambda T}U_T)^{-1}$. Hence the partial sums
of~\eqref{eq:CCC} are the rules $\Qt{2^{-N}}$, and~\eqref{eq:CCC} is their telescoped
limit. For $A=0$ the scalar identity behind~\eqref{eq:CCC},
\begin{equation}\label{eq:viete}
\frac1p=\frac{1}{1-e^{-p}}-\sum_{k\ge1}\frac{2^{-k}}{1+e^{-p2^{-k}}},
\end{equation}
is the logarithmic derivative of the hyperbolic form of Vi\`ete's product,
$\sinh x/x=\prod_{k\ge1}\cosh(x2^{-k})$, $x=p/2$ (cf.~\cite{Kac}), and the doubling step
is the case $s=0$ of the distribution relation of the polylogarithm~\cite{Lewin}.

The paper is organized around the following points; the first two are elementary
once the rectangle rule is identified, the third is the main result, and the last
translates the results for practitioners.

\begin{enumerate}
\item \emph{Generality} (\S\ref{sec:identity}). Nothing in the argument uses
self-adjointness: the identity holds for every $C_0$-semigroup on a Banach space,
with the dyadic ratio replaced by any sequence of integer ratios $b_k\ge2$
(Theorem~\ref{thm:mixed}). The only analytic input is the convergence of Riemann
sums.

\item \emph{Aliasing} (\S\ref{sec:aliasing}). The error of $\Qt{T}$ is the scalar
$T/2$ plus the resolvent summed over the alias points $\lambda-2\pi in/T$, $n\ne0$
(Theorem~\ref{thm:alias}, Corollary~\ref{cor:error}). This is the operator form of the
classical aliasing formula for impulse invariance, including the midpoint correction
of~\cite{Jackson}; for scalar integrands it is the Poisson summation formula behind
the error analysis of the trapezoidal rule~\cite{TW}. Specific to the dyadic
setting is that the $k$-th level removes exactly the aliases of $2$-adic valuation
$k-1$ (Corollary~\ref{cor:levels}). Identifiability statements such
as~\cite[Cor.~5]{CastilloV1} need no resolvent identity (Remark~\ref{rem:ident}).

\item \emph{Rates} (\S\ref{sec:rates}). For smooth data an Euler--Maclaurin
expansion (Theorem~\ref{thm:EM}) shows that the leading error is the known scalar
$T/2$ and justifies Romberg extrapolation~\cite{Romberg}
(Corollary~\ref{cor:romberg}). For data of limited regularity, pairing the aliases
writes the corrected error as a sum of resolvent differences on the line
$\Re z=\Re\lambda$ (Theorem~\ref{thm:aliasline}). The main result,
Theorem~\ref{thm:analytic}, gives the rates $T^\alpha$ in operator norm,
$T^{\alpha+\theta}$ on the interpolation spaces $X_\theta$, and $T^{1+\alpha}$ on
$D(B)$ for semigroups in the Crandall--Pazy class with parameter $\alpha\in(0,1]$,
i.e.\ $\|BS_t\|=O(t^{-1/\alpha})$; the last rate needs only resolvent decay of order
$\alpha$ along the line. Every exponent is attained
(Proposition~\ref{prop:analytic-sharp}), so the rates form a sharp scale from
unitary groups ($\alpha=0$, Theorem~\ref{thm:sa}) to analytic semigroups
($\alpha=1$); see Remark~\ref{rem:dichotomy}.

\item \emph{Spectral measures} (\S\ref{sec:spectral}). For self-adjoint $A$, the
contribution of each dyadic level to Stone's formula has a closed form in terms of
translates of the spectral measure (Theorem~\ref{thm:spectral}). The recovery of
$\int f\,d\mu_\psi$ itself is Stone's formula; the point is the accounting, which
shows in particular that the hypotheses of~\cite[Prop.~8, Cor.~4, Prop.~9]{CastilloV1}
are unnecessary.

\item \emph{Practice} (\S\ref{sec:practical}). For Green functions computed from a
sampled autocorrelation record, the results give an explicit error bound with a
Nyquist-type step-size rule (Proposition~\ref{prop:practical}), end corrections from
energy moments, Romberg extrapolation from a single record, and the record length
needed for a given broadening.
\end{enumerate}

\emph{Relation to prior work.} Rates of convergence for approximations of semigroups
by functions of their generators are classical: rational approximations and
time-stepping schemes~\cite{HershKato,BrennerThomee,Kovacs,BGT}, and approximation
formulas of Euler, Yosida and Dunford--Segal type~\cite{GomilkoTomilov,GKT}. A
recurring theme there is a dichotomy: for general bounded semigroups convergence
holds only on domains and interpolation spaces, while for holomorphic semigroups it
holds in operator norm and gains an order on every regularity class.
Theorem~\ref{thm:analytic} is the counterpart of this dichotomy for the sampled
resolvent $\Qt{T}$, a function of the generator whose symbol is the Laplace
transform of a signed measure. Those results concern approximations of the semigroup
itself, while $\Qt{T}$ approximates the resolvent, so they do not apply directly; in
the language of exponential integrators~\cite{HO},
$\Qt{T}=(\lambda-B)^{-1}g\bigl(T(\lambda-B)\bigr)$, where $g(z)=z/(1-e^{-z})$ is the
reciprocal of the first $\varphi$-function evaluated at $-z$.
Theorem~\ref{thm:analytic} also covers the intermediate Crandall--Pazy
class~\cite{CrandallPazy}. For a different discretization on that class, the decay
of Crank--Nicolson iterates, Wakaiki obtained rates interpolating between the
general and holomorphic cases~\cite{Wakaiki24}; on Hilbert spaces the class is
characterized by resolvent decay along a vertical line~\cite{Wakaiki24,Wakaiki25}.
In sampled-data control, analyticity plays a similar role in the robustness of
stabilization under fast sampling~\cite{LRT}.

As a quadrature rule for a Laplace integral whose integrand can be evaluated at
arbitrary times, the uniform rectangle rule is not competitive with exponentially
convergent rules~\cite{TW,WT}. The sampled resolvent is natural when the propagator
is available only on a uniform clock, and in that setting the alias accounting and the
rates of \S\ref{sec:rates} describe what can be achieved;
Section~\ref{sec:discussion} discusses this further. The Mellin-transform analogues
of the present results---fractional powers, spectral zeta functions and
determinants---will be treated separately.

Parts of the development in~\cite{CastilloV1} were carried out in the author's
doctoral dissertation~\cite{CastilloThesis}.

%% =====================================================================
\section{Setting: the sampled resolvent as a rectangle rule}\label{sec:setting}

Throughout, $X$ is a complex Banach space and $(S_t)_{t\ge0}$ a $C_0$-semigroup on
$X$ with generator $B$ and growth bound $\omega_0$. We fix $\lambda\in\C$ with
$\Re\lambda>\omega_0$, and numbers $\omega<\Re\lambda$ and $M\ge1$ with
$\|S_t\|\le Me^{\omega t}$ for $t\ge0$ (any $\omega>\omega_0$ admits such an $M$; for
unitary groups one may take $\omega=0$, $M=1$), and put $\delta:=\Re\lambda-\omega>0$. We use
the standard facts~\cite[Ch.~II, Ch.~IV]{EN}
\begin{equation}\label{eq:laplace}
(\lambda-B)^{-1}x=\int_0^\infty e^{-\lambda t}S_tx\,dt,\qquad
\|(\lambda-B)^{-1}\|\le \frac{M}{\delta},\qquad
r(S_T)=e^{\omega_0T},
\end{equation}
where $r$ denotes the spectral radius and the integral is an absolutely convergent
Bochner integral. We refer to these hypotheses as \textbf{(S)}.

\begin{definition}\label{def:FQ}
For $T>0$ and $\zeta\in\C$ with $|\zeta|=1$ put
\[
\F{T}{\zeta}:=\bigl(I-\zeta e^{-\lambda T}S_T\bigr)^{-1}
=\sum_{m=0}^{\infty}\zeta^me^{-\lambda mT}S_{mT},\qquad
\Qt{T}:=T\,\F{T}{1},\qquad
\Qtr{T}:=\Qt{T}-\tfrac T2\,I .
\]
\end{definition}

The inverse exists and the Neumann series converges in operator norm, because
$r(\zeta e^{-\lambda T}S_T)=e^{(\omega_0-\Re\lambda)T}<1$ by~\eqref{eq:laplace}.
With $h(t):=e^{-\lambda t}S_t$, $\Qt{T}=T\sum_{m\ge0}h(mT)$ is the left-endpoint
rectangle rule with step $T$ for~\eqref{eq:laplace}, and
$\Qtr{T}=T\bigl[\tfrac12h(0)+\sum_{m\ge1}h(mT)\bigr]$ is the trapezoid rule.

\begin{lemma}\label{lem:bound}
Under \textup{(S)}, $\|\Qt{T}\|\le M\,(T+1/\delta)$ for every $T>0$.
\end{lemma}

\begin{proof}
$\|\Qt{T}\|\le TM\sum_{m\ge0}e^{-\delta mT}=MT/(1-e^{-\delta T})\le M(T+1/\delta)$,
using $y/(1-e^{-y})\le1+y$ for $y>0$.
\end{proof}

\begin{proposition}[Riemann sums of the Laplace integral]\label{prop:riemann}
Assume \textup{(S)}.
\begin{enumerate}
\item[(i)] For $x\in D(B)$ and $T>0$,
$\ \|\Qt{T}x-(\lambda-B)^{-1}x\|\le \dfrac{MT}{\delta}\,\|(\lambda-B)x\|$.
\item[(ii)] For every $x\in X$, $\Qt{T}x\to(\lambda-B)^{-1}x$ as $T\to0^+$.
\item[(iii)] Let $0<\theta<1$ and $X_\theta:=(X,D(B))_{\theta,\infty}$ be the real
interpolation space, with $D(B)$ carrying the graph norm. There is $C$, depending
only on $M,\delta,|\lambda|$, such that
$\|\Qt{T}x-(\lambda-B)^{-1}x\|\le C\,T^{\theta}\|x\|_{X_\theta}$ for $x\in X_\theta$
and $0<T\le1$.
\end{enumerate}
\end{proposition}

\begin{proof}
(i) For $x\in D(B)$ the function $h(t)x=e^{-\lambda t}S_tx$ is continuously
differentiable with derivative $e^{-\lambda t}S_t(B-\lambda)x$, of norm at most
$Me^{-\delta t}\|(\lambda-B)x\|$. Therefore
\[
\Bigl\|\Qt{T}x-\int_0^\infty h(t)x\,dt\Bigr\|
\le\sum_{m\ge0}\int_{mT}^{(m+1)T}\!\!\int_{mT}^{t}\|(hx)'(s)\|\,ds\,dt
\le T\int_0^\infty\|(hx)'(s)\|\,ds\le\frac{MT}{\delta}\|(\lambda-B)x\| .
\]
(ii) follows from (i), the density of $D(B)$, and the uniform bound of
Lemma~\ref{lem:bound} for $T\le1$.
(iii) Write $E_T:=\Qt{T}-(\lambda-B)^{-1}$, so $\|E_T\|\le M(1+2/\delta)$ for $T\le1$
by Lemma~\ref{lem:bound} and~\eqref{eq:laplace}. For $y\in D(B)$, (i) gives
$\|E_Tx\|\le M(1+2/\delta)\|x-y\|+(MT/\delta)(1+|\lambda|)\|y\|_{D(B)}$. Taking the
infimum over $y$ bounds $\|E_Tx\|$ by a constant times the $K$-functional
$K(T,x)$, which is at most $T^\theta\|x\|_{X_\theta}$ by definition~\cite{Lunardi}.
\end{proof}

\begin{remark}[The scalar shadow]\label{rem:scalar}
For $X=\C$, $S_t\equiv1$ and $\lambda=p>0$ one has $\Qt{T}=T/(1-e^{-pT})$, and the
results of \S\ref{sec:identity} reduce to~\eqref{eq:viete}. Replacing the weight
$e^{-pt}$ by $t^{-s}e^{-pt}$, $s<1$, the rectangle rule with step $T$, with the
node $t=0$ omitted (it is singular for $0<s<1$), becomes $T^{1-s}\Lis_s(e^{-pT})$
and the same telescoping yields
$\Gamma(1-s)p^{s-1}=\Lis_s(e^{-p})-\sum_{k\ge1}2^{-k(1-s)}\Lis_s(-e^{-p2^{-k}})$,
which is~\cite[Lemma~21]{CCC} and~\cite[(37)]{CastilloV1}; the doubling step is the
distribution relation $\Lis_s(z)+\Lis_s(-z)=2^{1-s}\Lis_s(z^2)$~\cite{Lewin}. These
Mellin variants underlie the spectral zeta functions of~\cite{CastilloV1} and are
not pursued here.
\end{remark}

%% =====================================================================
\section{The distribution relation and the mixed-radix identity}\label{sec:identity}

\begin{lemma}[Distribution relation]\label{lem:dist}
Assume \textup{(S)}. Let $b\ge2$ be an integer, $\xi_b:=e^{2\pi i/b}$, $T>0$ and
$|\zeta|=1$. Then
\[
\F{bT}{\zeta^b}=\frac1b\sum_{r=0}^{b-1}\F{T}{\zeta\xi_b^r},
\qquad\text{in particular}\qquad
\Qt{bT}=\Qt{T}+T\sum_{r=1}^{b-1}\F{T}{\xi_b^r}.
\]
\end{lemma}

\begin{proof}
Expand each $\F{T}{\zeta\xi_b^r}$ in its norm-convergent Neumann series and use
$\frac1b\sum_{r=0}^{b-1}\xi_b^{rm}=1$ if $b\mid m$ and $0$ otherwise:
$\frac1b\sum_r\F{T}{\zeta\xi_b^r}
=\sum_{l\ge0}\zeta^{bl}e^{-\lambda lbT}S_{lbT}=\F{bT}{\zeta^b}$.
The second identity is the case $\zeta=1$ multiplied by $bT$.
\end{proof}

For $b=2$ and $\zeta=1$ the lemma reads $\F{2T}{1}=\frac12(\F{T}{1}+\F{T}{-1})
=\F{T}{1}\F{T}{-1}$, the ``Landen recursion'' of~\cite[Prop.~4]{CastilloV1}; it is
the operator form of the distribution relation of the geometric series, and is
unrelated to Landen's transformation of elliptic integrals.

\begin{theorem}[Mixed-radix resolvent identity]\label{thm:mixed}
Assume \textup{(S)}. Let $T_0>0$, let $(b_k)_{k\ge1}$ be integers $\ge2$, and put
$T_k:=T_{k-1}/b_k$ and $\xi_k:=e^{2\pi i/b_k}$. Then for every $N\ge0$
\begin{equation}\label{eq:partial}
\Qt{T_0}-\sum_{k=1}^{N}T_k\sum_{r=1}^{b_k-1}\F{T_k}{\xi_k^r}=\Qt{T_N},
\end{equation}
and consequently
\begin{equation}\label{eq:mixed}
(\lambda-B)^{-1}=\Qt{T_0}-\sum_{k=1}^{\infty}T_k\sum_{r=1}^{b_k-1}\F{T_k}{\xi_k^r}
\end{equation}
in the strong operator topology. The partial sums obey the estimates of
Proposition~\ref{prop:riemann} with $T=T_N$. If $B$ is bounded, \eqref{eq:mixed}
converges in operator norm and $\|\Qt{T_N}-(\lambda-B)^{-1}\|\le(MT_N/\delta)\|\lambda-B\|$.
\end{theorem}

\begin{proof}
Lemma~\ref{lem:dist} with $T=T_k$ and $b=b_k$ gives
$\Qt{T_{k-1}}-\Qt{T_k}=T_k\sum_{r=1}^{b_k-1}\F{T_k}{\xi_k^r}$; summing over
$k=1,\dots,N$ gives~\eqref{eq:partial}. Since $T_N\le2^{-N}T_0\to0$, the remaining
assertions are Proposition~\ref{prop:riemann}, applied with
$\|(\lambda-B)x\|\le\|\lambda-B\|\,\|x\|$ in the bounded case.
\end{proof}

\begin{corollary}[{\cite[Prop.~22(i)]{CCC}}]\label{cor:CCC}
Let $A$ be self-adjoint on a Hilbert space and $\lambda>0$. Then~\eqref{eq:CCC}
holds in the strong operator topology.
\end{corollary}

\begin{proof}
Take $S_t=U_t=e^{-itA}$, so $B=-iA$, $\omega_0=0$, and $\omega=0$, $M=1$, $\delta=\lambda$, and note
$(\lambda-B)^{-1}=(\lambda+iA)^{-1}=-i(A-i\lambda)^{-1}$. Apply
Theorem~\ref{thm:mixed} with $T_0=1$ and $b_k\equiv2$, so that $\xi_k=-1$, and
multiply by $i$.
\end{proof}

The only analytic input is the convergence of Riemann sums,
Proposition~\ref{prop:riemann}; the argument of~\cite{CCC} proceeds instead
through a dyadic expansion of the Cauchy kernel and the spectral theorem. The
derivation of~\cite[Thm.~1]{CastilloV1} is the case $b_k\equiv2$ of the present one
specialized to unitary groups.

\begin{remark}[Order of summation]\label{rem:order}
Expanding each $\F{T_k}{\xi_k^r}$ in its Neumann series gives the double-series
form used in~\cite{CCC,CastilloV1}. The levels are not absolutely summable in
general: for $B=-iA$ with $\sigma(A)=\R$ and $b_k\equiv2$,
$\|T_k\F{T_k}{-1}\|=T_k\sup_{a\in\R}|1+e^{-T_k(\lambda+ia)}|^{-1}
=T_k/(1-e^{-\lambda T_k})\ge1/\lambda$, so $\sum_k\|T_k\F{T_k}{-1}\|=\infty$. Only the
telescoped partial sums~\eqref{eq:partial} converge, which is why the limits
in~\cite{CCC} cannot be interchanged.
\end{remark}

\begin{remark}[Evaluation]\label{rem:stability}
For bounded $B$, computing $\Qt{T_N}$ directly requires inverting
$I-e^{-\lambda T_N}S_{T_N}$, whose norm is $O(T_N)$, so the subtraction loses accuracy
as $T_N\to0$. In the telescoped
form~\eqref{eq:partial} the only inversion with $\zeta=1$ occurs at the coarse step
$T_0$; for $k\ge1$ and bounded $B$ one has
$I-\zeta e^{-\lambda T_k}S_{T_k}\to(1-\zeta)I$ with $\zeta\ne1$, which is well
conditioned. The effect is modest in practice (Section~\ref{sec:nonnormal}). When the
samples are produced by scaling and squaring~\cite{Higham}, the coarse samples
$S_{T_k}$ are by-products of the computation of the finest one.
\end{remark}

%% =====================================================================
\section{Aliasing}\label{sec:aliasing}

\begin{theorem}[Poisson summation for the sampled resolvent]\label{thm:alias}
Assume \textup{(S)}. Let $T>0$, $\eta\in\R$ and $\zeta=e^{i\eta}$. For every
$x\in D(B)$,
\begin{equation}\label{eq:alias}
T\F{T}{\zeta}x=\frac T2\,x
+\lim_{L\to\infty}\sum_{|n|\le L}\Bigl(\lambda-\frac{i(\eta+2\pi n)}{T}-B\Bigr)^{-1}x .
\end{equation}
If $B$ is bounded, \eqref{eq:alias} holds in operator norm.
\end{theorem}

Every resolvent on the right exists, since
$\Re\bigl(\lambda-i(\eta+2\pi n)/T\bigr)=\Re\lambda>\omega_0$.

\begin{proof}
Put $\lambda_\eta:=\lambda-i\eta/T$ and $g(t):=e^{-\lambda_\eta t}S_tx$ for
$t\ge0$, $g(t):=0$ for $t<0$. Since $\zeta^me^{-\lambda mT}=e^{-\lambda_\eta mT}$,
we have $T\F{T}{\zeta}x=T\sum_{m\ge0}g(mT)$. The periodization
$P(t):=\sum_{m\in\Z}g(t+mT)$ converges absolutely ($\|g(t)\|\le Me^{-\delta t}\|x\|$),
is $T$-periodic, and is continuous except at $T\Z$, where
$P(0^+)=\sum_{m\ge0}g(mT)$ and $P(0^-)=\sum_{m\ge1}g(mT)$. Its Fourier
coefficients are
\[
c_n=\frac1T\int_0^TP(t)e^{-2\pi int/T}dt=\frac1T\int_0^\infty g(t)e^{-2\pi int/T}dt
=\frac1T\Bigl(\lambda_\eta+\frac{2\pi in}{T}-B\Bigr)^{-1}x .
\]
Because $x\in D(B)$, $g$ is Lipschitz on $[0,\infty)$ with
$\|g'(t)\|\le Me^{-\delta t}\|(\lambda_\eta-B)x\|$, whence
$\|P(t)-P(0^+)\|\le Ct$ and $\|P(-t)-P(0^-)\|\le Ct$ for $0<t<T$, with
$C=M\|(\lambda_\eta-B)x\|/(1-e^{-\delta T})$. With the Dirichlet kernel,
\[
\sum_{|n|\le L}c_n-\frac{P(0^+)+P(0^-)}2
=\frac1T\int_0^{T/2}\frac{P(t)+P(-t)-P(0^+)-P(0^-)}{\sin(\pi t/T)}
\,\sin\!\Bigl(\frac{(2L+1)\pi t}{T}\Bigr)dt ,
\]
and the quotient is a bounded, hence Bochner integrable, $X$-valued function. The
Riemann--Lebesgue lemma holds for Bochner integrable functions (approximate by
$X$-valued step functions), so the right side tends to $0$ as $L\to\infty$. Thus
$\sum_{|n|\le L}c_n\to\sum_{m\ge0}g(mT)-\frac12x$; multiplying by $T$ and replacing
$n$ by $-n$ gives~\eqref{eq:alias}. If $B$ is bounded, the same argument applies to
the $\mathcal L(X)$-valued function $t\mapsto e^{-\lambda_\eta t}S_t$, which is
Lipschitz in operator norm.
\end{proof}

\begin{remark}[Prior forms]\label{rem:impulse}
For $X=\C$ and $S_t=e^{bt}$, \eqref{eq:alias} with $\eta=0$ is the aliasing
formula for the impulse-invariant discretization of the transfer function
$1/(\lambda-b)$~\cite{OS}, and the term $T/2$ is the correction of~\cite{Jackson} for
impulse responses that jump at $t=0$. For general scalar integrands it is the Poisson
summation formula underlying the error analysis of the trapezoidal rule~\cite{TW}.
Theorem~\ref{thm:alias} is the operator version, and its proof is the standard one;
what is used below is its form along the line $\Re z=\Re\lambda$
(Theorem~\ref{thm:aliasline}) and the level structure of
Corollary~\ref{cor:levels}.
\end{remark}

\begin{corollary}[Exact error representation]\label{cor:error}
Assume \textup{(S)}. For $x\in D(B)$ and $T>0$,
\begin{equation}\label{eq:error}
\Qt{T}x-(\lambda-B)^{-1}x
=\frac T2\,x+\lim_{L\to\infty}\sum_{0<|n|\le L}
\Bigl(\lambda-\frac{2\pi in}{T}-B\Bigr)^{-1}x .
\end{equation}
In particular, the truncation error of~\eqref{eq:mixed} after $N$ levels consists of
the scalar term $\frac{T_N}{2}x$ and the resolvent summed over the alias frequencies
$2\pi n/T_N$, $n\ne0$.
\end{corollary}

\begin{proof}
Take $\eta=0$ in Theorem~\ref{thm:alias} and isolate $n=0$; the second sentence
is~\eqref{eq:partial}.
\end{proof}

For a normal generator the alias terms are small when the spectral data live well
inside the Nyquist band $|\Im(\lambda-b)|<\pi/T$, and the scalar $T/2$ dominates the
error; this is quantified in \S\ref{sec:rates}.

The levels of~\eqref{eq:mixed} admit the same description. Let $\nu_0:=2\pi/T_0$,
$\beta_0:=1$, $\beta_k:=b_1\cdots b_k$, and
\[
\Lambda_k:=\{\beta_{k-1}m:\ m\in\Z,\ b_k\nmid m\}\qquad(k\ge1).
\]
The sets $\Lambda_k$ partition $\Z\setminus\{0\}$: for $\ell\ne0$ there is a largest
$j\ge0$ with $\beta_j\mid \ell$, and $\ell\in\Lambda_{j+1}$. For $b_k\equiv b$, $\ell\in\Lambda_k$
means that the $b$-adic valuation of $\ell$ is $k-1$.

\begin{corollary}[Level-by-level alias removal]\label{cor:levels}
Assume \textup{(S)} and let $x\in D(B)$. In the notation of Theorem~\ref{thm:mixed},
\begin{align*}
\Qt{T_0}x&=\frac{T_0}{2}x+\lim_{L\to\infty}\sum_{|\ell|\le L}(\lambda-i\nu_0\ell-B)^{-1}x,\\
T_k\sum_{r=1}^{b_k-1}\F{T_k}{\xi_k^r}x
&=\frac{T_{k-1}-T_k}{2}x+\sum_{r=1}^{b_k-1}\lim_{L\to\infty}\sum_{|n|\le L}
\bigl(\lambda-i\nu_0\beta_{k-1}(r+b_kn)-B\bigr)^{-1}x \qquad(k\ge1),
\end{align*}
and the frequencies $\nu_0\beta_{k-1}(r+b_kn)$, $1\le r<b_k$, $n\in\Z$, are exactly
$\nu_0\Lambda_k$.
\end{corollary}

\begin{proof}
Apply Theorem~\ref{thm:alias} with $T=T_k$ and $\eta=2\pi r/b_k$, and note
$(\eta+2\pi n)/T_k=2\pi(r+b_kn)/T_{k-1}=\nu_0\beta_{k-1}(r+b_kn)$ and
$(b_k-1)T_k=T_{k-1}-T_k$.
\end{proof}

\begin{remark}[Reading the identity]\label{rem:reading}
Formally, subtracting the level expansions from that of $\Qt{T_0}$, the constants
cancel because $\sum_{k\ge1}(T_{k-1}-T_k)=T_0$, and the aliases cancel because the
$\Lambda_k$ partition $\Z\setminus\{0\}$; only the term $\ell=0$, the resolvent,
survives. This formal computation rearranges conditionally convergent sums;
Theorem~\ref{thm:mixed} is the rigorous version, and Corollary~\ref{cor:error}
identifies precisely what remains after finitely many levels. The ``octave''
structure of~\eqref{eq:CCC} is thus the $2$-adic valuation of the alias index.
\end{remark}

\begin{remark}[Zak transforms and filter banks]\label{rem:zak}
In~\cite[Prop.~2]{CastilloV1} the matrix coefficients of the levels were read as
values $Z_TG(0,w)$ of Zak transforms of
$G(t)=i\ip{U_t\varphi}{\psi}e^{-\lambda t}e^{2\pi it}\mathbf 1_{[0,\infty)}(t)$.
The second Zak variable enters only through $\zeta=e^{2\pi i(T-w)}$, so by
Theorem~\ref{thm:alias} it sweeps the real part of the spectral parameter modulo
$2\pi/T$: at sampling time $T$ the Zak transform is the aliased resolvent. The value
at $t=0$, where $G$ jumps, is a convention; the right-continuous choice
of~\cite{CastilloV1} produces $\Qt{T}$, while the midpoint choice produces the
trapezoid rule $\Qtr{T}$, which is second order (\S\ref{sec:rates}). In the
filter-bank language of~\cite[Rem.~8]{CastilloV1}, level $k$ is a bank of $b_k-1$
one-pole filters whose passbands are the alias frequencies $\nu_0\Lambda_k$.
\end{remark}

\begin{remark}[Identifiability]\label{rem:ident}
If two $C_0$-semigroups satisfy $S^A_{T_k}=S^B_{T_k}$ along a sequence $T_k\to0$,
then they agree at every time $mT_k$, a dense subset of $[0,\infty)$, hence
everywhere by strong continuity, and their generators coincide. This
gives~\cite[Cor.~5]{CastilloV1} without the resolvent identity, for any null
sequence of times. Moreover $(S_{T_k})_{k\le N}$ is determined by $S_{T_N}$ alone when
the ratios are integers. The content of~\eqref{eq:mixed} is therefore not
identifiability but the organization of the truncation error described by
Corollaries~\ref{cor:error} and~\ref{cor:levels}.
\end{remark}

%% =====================================================================
\section{Rates and acceleration}\label{sec:rates}

By~\eqref{eq:partial}, every statement about the partial sums of~\eqref{eq:mixed} is
a statement about the rectangle rule $\Qt{T_N}$, so the classical theory of the
trapezoid rule applies. We write $B_{2j}$ for the Bernoulli numbers
($B_2=\frac16$, $B_4=-\frac1{30}$, \dots).

\begin{theorem}[Euler--Maclaurin expansion]\label{thm:EM}
Assume \textup{(S)}. Let $q\ge1$ and $x\in D(B^{2q})$. Then for every $T>0$
\begin{equation}\label{eq:EM}
\Qtr{T}x=(\lambda-B)^{-1}x+\sum_{j=1}^{q-1}\frac{B_{2j}}{(2j)!}\,T^{2j}(\lambda-B)^{2j-1}x
+\rho_q(T)x,
\qquad
\|\rho_q(T)x\|\le\frac{2|B_{2q}|}{(2q)!}\,\frac{M}{\delta}\,T^{2q}\,\|(\lambda-B)^{2q}x\| .
\end{equation}
\end{theorem}

\begin{proof}
Let $\phi(j):=T\,h(jT)x=Te^{-\lambda jT}S_{jT}x$. This is a $C^{2q}$ function of
$j\ge0$, and $\phi^{(i)}(j)=T^{i+1}e^{-\lambda jT}S_{jT}(B-\lambda)^ix$. The Euler--Maclaurin formula
with remainder~\cite[\S2.10]{DLMF},
\[
\sum_{j=0}^{n}\phi(j)=\int_0^n\phi+\frac{\phi(0)+\phi(n)}2
+\sum_{s=1}^{q-1}\frac{B_{2s}}{(2s)!}\bigl(\phi^{(2s-1)}(n)-\phi^{(2s-1)}(0)\bigr)
+\int_0^n\frac{B_{2q}-\widetilde B_{2q}(u)}{(2q)!}\,\phi^{(2q)}(u)\,du,
\]
with $\widetilde B_{2q}$ the periodic Bernoulli function, holds for $X$-valued $\phi$:
apply it to $x^*\circ\phi$ for each $x^*\in X^*$ and use Hahn--Banach. Let
$n\to\infty$; all boundary terms at $n$ vanish by the exponential decay, the
integral tends to $\int_0^\infty h(t)x\,dt=(\lambda-B)^{-1}x$, and
$-\phi^{(2s-1)}(0)=T^{2s}(\lambda-B)^{2s-1}x$. Since
$|B_{2q}-\widetilde B_{2q}(u)|\le2|B_{2q}|$, the remainder is bounded by
$\frac{2|B_{2q}|}{(2q)!}\,T^{2q}\int_0^\infty\|e^{-\lambda t}S_t(B-\lambda)^{2q}x\|\,dt$.
Finally $\sum_{j\ge0}\phi(j)-\frac12\phi(0)=\Qtr{T}x$.
\end{proof}

In particular $\Qt{T}x=(\lambda-B)^{-1}x+\frac T2x+O(T^2)$ for $x\in D(B^2)$: the
leading error of the dyadic identity is the scalar $\frac{T_N}{2}x$, which is known
and can be removed. For a scalar generator the expansion is the Taylor series of
$z/(1-e^{-z})=1+\frac z2+\sum_{j\ge1}B_{2j}z^{2j}/(2j)!$ at $z=T(\lambda-b)$; its
radius of convergence $2\pi$ is set by the poles $z=2\pi in$, which are the aliases
of Corollary~\ref{cor:error}.

\begin{corollary}[Trapezoid correction and Romberg extrapolation]\label{cor:romberg}
Assume \textup{(S)} and let $b_k\equiv b\ge2$, so $T_N=b^{-N}T_0$. Put
$P_N^{(0)}:=\Qtr{T_N}$, i.e.\ the partial sum~\eqref{eq:partial} minus
$\frac{T_N}{2}I$, and for $1\le j\le N$
\[
P_N^{(j)}:=\frac{b^{2j}P_N^{(j-1)}-P_{N-1}^{(j-1)}}{b^{2j}-1}.
\]
Then for $j\ge0$ and $x\in D(B^{2j+2})$,
$\|P_N^{(j)}x-(\lambda-B)^{-1}x\|\le C_{j,b}\,(M/\delta)\,T_N^{2j+2}\,
\|(\lambda-B)^{2j+2}x\|$, with $C_{j,b}$ depending only on $j$ and $b$.
\end{corollary}

\begin{proof}
By induction, $P_N^{(j)}=\sum_{i=0}^{j}\alpha_iP^{(0)}_{N-i}$ with
$\sum_i\alpha_i=1$ and $\sum_i\alpha_ib^{2li}=0$ for $1\le l\le j$. Apply
Theorem~\ref{thm:EM} with $q=j+1$ to each $P^{(0)}_{N-i}x$, noting
$T_{N-i}=b^iT_N$: the terms $T^{2l}$, $l\le j$, cancel, and
$P^{(j)}_Nx-(\lambda-B)^{-1}x=\sum_i\alpha_i\rho_{j+1}(b^iT_N)x$. One may take
$C_{j,b}=\frac{2|B_{2j+2}|}{(2j+2)!}\sum_i|\alpha_i|b^{(2j+2)i}$.
\end{proof}

For a bounded generator Corollary~\ref{cor:romberg} holds in operator norm, with
$\|(\lambda-B)^{2j+2}\|$ in place of $\|(\lambda-B)^{2j+2}x\|$; the asymptotic regime
begins once $T_N(|\lambda|+\|B\|)$ is of order one. For unbounded generators the
error constants involve high powers of $B$, and the relevant question is how fast
$\Qt{T}$ converges on vectors of limited regularity. For self-adjoint generators
this has a sharp answer.

\begin{theorem}[Self-adjoint generators]\label{thm:sa}
Let $A$ be self-adjoint on a Hilbert space $\Hil$, $S_t=e^{-itA}$, $\lambda>0$ and
$0<T\le T_*$. Put $g(z):=z/(1-e^{-z})$,
$K:=\{z:0\le\Re z\le\lambda T_*,\ |\Im z|\le\pi\}$,
\[
c_1:=\sup_{z\in K}\Bigl|\frac{g(z)-1}{z}\Bigr|,\qquad
c_2:=\sup_{z\in K}\Bigl|\frac{g(z)-1-z/2}{z^2}\Bigr|,\qquad
c_3:=\frac2\lambda+T_*,
\]
which are finite because $g$ is analytic in $|\Im z|<2\pi$.
\begin{enumerate}
\item[(i)] If $\psi\in D(|A|^s)$ with $0<s\le1$, then
$\|(\Qt{T}-(\lambda+iA)^{-1})\psi\|\le c_1T\|\psi\|+c_3(T/\pi)^s\||A|^s\psi\|$.
\item[(ii)] If $\psi\in D(|A|^s)$ with $1\le s\le2$, then
$\|(\Qtr{T}-(\lambda+iA)^{-1})\psi\|\le c_2T^2\|(\lambda+iA)\psi\|
+(c_3+T_*/2)(T/\pi)^s\||A|^s\psi\|$.
\item[(iii)] The exponent $s$ in \textup{(i)} and \textup{(ii)} cannot be improved:
for every $s'>0$ there are $A$ and $\psi\in\bigcap_{s<s'}D(|A|^s)$ with
$\|(\Qt{T_N}-(\lambda+iA)^{-1})\psi\|\ge c\,T_N^{s'}$ and
$\|(\Qtr{T_N}-(\lambda+iA)^{-1})\psi\|\ge c\,T_N^{s'}$ for $T_N=2^{-N}$ and all large
$N$, with $c>0$.
\end{enumerate}
\end{theorem}

\begin{proof}
By the spectral theorem, $\Qt{T}-(\lambda+iA)^{-1}$ is multiplication by
$e_T(a):=T/(1-e^{-T(\lambda+ia)})-1/(\lambda+ia)$ in a spectral representation of
$A$, and $\|(\Qt{T}-(\lambda+iA)^{-1})\psi\|^2=\int|e_T|^2d\mu_\psi$. Write
$p=\lambda+ia$ and $z=Tp$. On the band $|a|\le\pi/T$ one has $z\in K$ and
$e_T(a)=T(g(z)-1)/z$, so $|e_T(a)|\le c_1T$; likewise
$|e_T(a)-T/2|=T^2|p|\,|(g(z)-1-z/2)/z^2|\le c_2T^2|p|$. Off the band,
$|1-e^{-z}|\ge1-e^{-\lambda T}$ gives
$|e_T(a)|\le T/(1-e^{-\lambda T})+1/\lambda\le c_3$. By Chebyshev's inequality,
$\mu_\psi(|a|>\pi/T)\le(T/\pi)^{2s}\||A|^s\psi\|^2$. Combining the two bands and
using $\sqrt{u^2+v^2}\le u+v$ gives (i) and (ii).

For (iii), let $A$ have simple eigenvalues $a_j=2\pi2^j$, $j\ge0$, with orthonormal
eigenvectors $e_j$, and $\psi=\sum_j2^{-js'}e_j$; then $\psi\in D(|A|^s)$ iff $s<s'$.
For $j\ge N$, $T_Na_j\in2\pi\Z$, so
$|e_{T_N}(a_j)|\ge T_N/(1-e^{-\lambda T_N})-1/a_j\ge\frac1\lambda-\frac1{2\pi 2^N}$,
which is at least $\frac1{2\lambda}$ for large $N$, and the same holds for
$e_{T_N}(a_j)-T_N/2$ once $T_N\le\frac1{2\lambda}$ as well. Hence both errors are at
least $\frac1{4\lambda}\bigl(\sum_{j\ge N}2^{-2js'}\bigr)^{1/2}\ge c\,2^{-Ns'}$.
\end{proof}

The worst case in (iii) places spectral mass exactly on the aliases of every level,
where the symbol of $\Qt{T}$ has height $\approx1/\lambda$. This is the regime in which
the dyadic levels are genuinely needed: each level detects one octave of the
spectral tail.

\begin{remark}[Smooth spectral measures]\label{rem:smooth}
Away from the aliases the symbol of $\Qt{T}$ is small: its mean square over a period
$|a-2\pi n/T|\le\pi/T$ is $T^2/(1-e^{-2\lambda T})\approx T/(2\lambda)$. If
$\mu_\psi$ has a density that varies little on intervals of length $2\pi/T$ beyond
$|a|=\pi/T$, the off-band contribution is therefore of order
$(T\,\mu_\psi(|a|>\pi/T))^{1/2}$ rather than $\mu_\psi(|a|>\pi/T)^{1/2}$, which
suggests the rates $T^{\min(1,s+1/2)}$ for $\Qt{T}$ and $T^{\min(2,s+1/2)}$ for
$\Qtr{T}$. Figure~\ref{fig:rates}(b) is consistent with this; we do not pursue a
precise statement.
\end{remark}

\subsection{Non-normal generators: resolvent decay along the alias line}
\label{sec:nonnormal-rates}

Theorem~\ref{thm:sa} rests on the spectral theorem. For a general generator, the
right object is not the location of the spectrum but the size of the resolvent at
the points where the aliases sit. Pairing the terms $n$ and $-n$ in
Corollary~\ref{cor:error} turns the conditionally convergent alias sum into a sum
over the vertical line $\Re z=\Re\lambda$ that gains one power of the alias
frequency.

\begin{theorem}[Alias-line representation of the trapezoid error]\label{thm:aliasline}
Assume \textup{(S)}. For $x\in D(B)$ and $T>0$, with $\nu_n:=2\pi n/T$,
\begin{align}
\Qtr{T}x-(\lambda-B)^{-1}x
&=2\lim_{L\to\infty}\sum_{n=1}^{L}(\lambda-i\nu_n-B)^{-1}(\lambda+i\nu_n-B)^{-1}(\lambda-B)x
\label{eq:aliasline}\\
&=\lim_{L\to\infty}\sum_{n=1}^{L}\frac{1}{i\nu_n}
\Bigl[(\lambda-i\nu_n-B)^{-1}-(\lambda+i\nu_n-B)^{-1}\Bigr](\lambda-B)x .
\label{eq:aliasline2}
\end{align}
In particular, whenever the right-hand side is finite,
\begin{equation}\label{eq:aliasbound}
\|\Qtr{T}x-(\lambda-B)^{-1}x\|\le\sum_{n\ge1}\frac{1}{\nu_n}
\Bigl(\|(\lambda-i\nu_n-B)^{-1}\|+\|(\lambda+i\nu_n-B)^{-1}\|\Bigr)\|(\lambda-B)x\| .
\end{equation}
\end{theorem}

\begin{proof}
Write $R(\mu)=(\mu-B)^{-1}$. For $\mu,\mu'\in\rho(B)$,
$R(\mu)+R(\mu')=R(\mu)\bigl[(\mu'-B)+(\mu-B)\bigr]R(\mu')$ and
$R(\mu)-R(\mu')=(\mu'-\mu)R(\mu)R(\mu')$. With $\mu=\lambda-i\nu$,
$\mu'=\lambda+i\nu$ the bracket is $2(\lambda-B)$ and $\mu'-\mu=2i\nu$; for
$x\in D(B)$ the operators $\lambda-B$ and $R(\mu')$ commute on $x$. Hence
$R(\lambda-i\nu_n)x+R(\lambda+i\nu_n)x=2R(\lambda-i\nu_n)R(\lambda+i\nu_n)(\lambda-B)x
=(i\nu_n)^{-1}[R(\lambda-i\nu_n)-R(\lambda+i\nu_n)](\lambda-B)x$, and
\eqref{eq:aliasline}--\eqref{eq:aliasline2} are~\eqref{eq:error} with the symmetric
partial sums grouped in pairs.
\end{proof}

For a scalar generator $b$, with $p=\lambda-b$, the identity is the partial-fraction
expansion
\[
\sum_{n\ge1}\frac{2p}{p^2+\nu_n^2}=\frac T2\coth\frac{pT}2-\frac1p .
\]
For a unitary group $e^{-itA}$ with $\sigma(A)=\R$, $\|(\lambda\pm i\nu-B)^{-1}\|=1/\lambda$
for every $\nu$, the series in~\eqref{eq:aliasbound} diverges logarithmically, and
the bound is void, in accordance with Theorem~\ref{thm:sa}(iii). Any decay of the
resolvent along the line $\Re z=\Re\lambda$ makes it summable. We measure that decay,
and the corresponding smoothing in time, on two scales.

\begin{definition}\label{def:scales}
Let $0<\alpha\le1$. Condition \textup{(R$_\alpha$)} holds if there is $K_R\ge1$ with
\[
\|(\lambda-i\nu-B)^{-1}\|\le K_R|\nu|^{-\alpha}\qquad(|\nu|\ge1).
\]
Condition \textup{(D$_\alpha$)} holds if $S_tX\subset D(B)$ for $t>0$ and there are
$K\ge1$ and $\omega_1<\Re\lambda$ with
\[
\|(\lambda-B)S_t\|\le K\bigl(1+t^{-1/\alpha}\bigr)e^{\omega_1t}\qquad(t>0).
\]
We then put $\delta_1:=\Re\lambda-\omega_1>0$.
\end{definition}

Condition \textup{(D$_\alpha$)} says that $(S_t)$ belongs to the Crandall--Pazy
class with parameter $\alpha$~\cite{CrandallPazy}, in the terminology
of~\cite{Wakaiki24}: $S_tX\subset D(B)$ for $t>0$ and $\|BS_t\|=O(t^{-1/\alpha})$ as
$t\to0^+$. The two formulations differ by $\lambda S_t$, which is bounded for
$0<t\le1$. The case $\alpha=1$ is the class of analytic semigroups
(Remark~\ref{rem:oneline}). On Hilbert spaces \textup{(D$_\alpha$)} and
\textup{(R$_\alpha$)} are equivalent~\cite[Thm.~2.3]{Wakaiki24}. Since
$(\lambda-B)S_t=(\lambda-B)S_1S_{t-1}$ for $t\ge1$, it suffices in
\textup{(D$_\alpha$)} to have the bound for $0<t\le1$. A normal operator whose
spectrum lies in a region $\{z:\Re z\le\omega-c|\Im z|^\alpha\}$ satisfies both
conditions. Non-analytic examples with $0<\alpha<1$ arise from elastic systems
$\ddot u+\rho A^\beta\dot u+Au=0$ with $A>0$ self-adjoint and $0<\beta<\frac12$,
whose semigroups are of Gevrey class of every order greater than
$1/(2\beta)$~\cite{ChenTriggiani} (their $\alpha$ is our $\beta$); for
the model case of \S\ref{sec:num-damped} both conditions hold with $\alpha=2\beta$.

The implication \textup{(D$_\alpha$)}$\Rightarrow$\textup{(R$_\alpha$)} is known
(cf.~\cite[Thm.~4.1(a)]{Wakaiki25}); we include the short proof to record the constant.

\begin{lemma}\label{lem:DtoR}
Assume \textup{(S)} and \textup{(D$_\alpha$)}. Then \textup{(R$_\alpha$)} holds with
$K_R=2M+2K+K^2(2/\delta_1+c_\alpha)$, where $c_\alpha:=2^{1+2/\alpha}/(2/\alpha-1)$.
\end{lemma}

\begin{proof}
Let $h(t)=e^{-\lambda t}S_t$ and $x\in X$. By \textup{(D$_\alpha$)}, $S_tx\in D(B^2)$
for $t>0$ ($S_tx=S_{t/2}S_{t/2}x$ and $BS_tx=S_{t/2}BS_{t/2}x$), and $t\mapsto h(t)x$
is $C^2$ on $(0,\infty)$ with $h'(t)x=-(\lambda-B)h(t)x$ and
$h''(t)x=(\lambda-B)^2h(t)x$. Moreover $\|h(t)\|\le M$,
$\|h'(t)\|\le K(1+t^{-1/\alpha})e^{-\delta_1t}$ and, applying \textup{(D$_\alpha$)}
to $S_{t/2}$ twice,
$\|h''(t)\|\le K^2(1+(2/t)^{1/\alpha})^2e^{-\delta_1t}$, whence for every $\tau>0$
\begin{equation}\label{eq:h2int}
\int_\tau^\infty\|h''(t)\|\,dt\le2K^2\!\int_\tau^\infty\bigl(1+(2/t)^{2/\alpha}\bigr)e^{-\delta_1t}dt
\le K^2\bigl(2/\delta_1+c_\alpha\tau^{1-2/\alpha}\bigr).
\end{equation}
Let $|\nu|\ge1$ and $\tau:=|\nu|^{-\alpha}\le1$. In
$(\lambda-i\nu-B)^{-1}x=\int_0^\infty e^{i\nu t}h(t)x\,dt$ the integral over $[0,\tau]$
has norm at most $M\tau\|x\|$. On $[\tau,\infty)$ integrate by parts twice:
\[
\int_\tau^\infty e^{i\nu t}h(t)x\,dt
=-\frac{e^{i\nu\tau}h(\tau)x}{i\nu}+\frac{e^{i\nu\tau}h'(\tau)x}{(i\nu)^2}
+\frac1{(i\nu)^2}\int_\tau^\infty e^{i\nu t}h''(t)x\,dt ,
\]
the boundary terms at infinity vanishing by the exponential decay. Hence
\[
\|(\lambda-i\nu-B)^{-1}\|\le M\tau+\frac M{|\nu|}
+\frac{K(1+\tau^{-1/\alpha})+K^2(2/\delta_1+c_\alpha\tau^{1-2/\alpha})}{\nu^2}.
\]
With $\tau^{-1/\alpha}=|\nu|$ and $\tau^{1-2/\alpha}=|\nu|^{2-\alpha}$, every term is
at most its coefficient times $|\nu|^{-\alpha}$, because $|\nu|\ge1$ and $\alpha\le1$.
\end{proof}

\begin{remark}[Analyticity from one line]\label{rem:oneline}
For $\alpha=1$ the two conditions are equivalent, and both characterize analyticity:
under \textup{(S)}, $(S_t)$ is analytic if and only if \textup{(D$_1$)} holds, if and
only if \textup{(R$_1$)} holds. This is known, and we only indicate how it follows
from the literature. Analyticity and \textup{(D$_1$)} are equivalent by the standard
characterizations~\cite[Ch.~2, Thm.~5.2]{Pazy},~\cite[Thm.~II.4.6]{EN}, applied to
$e^{-\omega't}S_t$ with $\omega_0<\omega'<\Re\lambda$, and
\textup{(D$_1$)}$\Rightarrow$\textup{(R$_1$)} is Lemma~\ref{lem:DtoR}. Conversely,
assume \textup{(R$_1$)} and put $a':=\Re\lambda-\delta/(2M)$, so that
$a'\ge(\Re\lambda+\omega)/2>\omega$. Since $\|(\Re\lambda+is-B)^{-1}\|\le M/\delta$, a
Neumann series gives
$\|(a'+is-B)^{-1}\|\le2\|(\Re\lambda+is-B)^{-1}\|$ for all $s\in\R$. Hence the
generator $B-a'$ of the bounded semigroup $e^{-a't}S_t$ satisfies
$\|(is-(B-a'))^{-1}\|\le C/|s|$ for all $s\ne0$, and by~\cite[Thm.~3.12]{HytonenLappas},
a variant of~\cite[Cor.~3.7.12]{ABHN} and~\cite[Ch.~2, Thm.~5.2]{Pazy}, $e^{-a't}S_t$ is a
bounded holomorphic semigroup. In particular
$\sup_{t>0}\|t(B-a')e^{-a't}S_t\|<\infty$, which gives \textup{(D$_1$)} with
$\omega_1=a'$. The point for us is that a resolvent bound on the single line
$\Re z=\Re\lambda$, where the aliases sit, already suffices.
\end{remark}

\begin{theorem}[Rates on the Crandall--Pazy scale]\label{thm:analytic}
Assume \textup{(S)}, let $0<\alpha\le1$ and $0<T\le1$.
\begin{enumerate}
\item[(i)] Under \textup{(D$_\alpha$)},
$\|\Qtr{T}-(\lambda-B)^{-1}\|\le C_1T^\alpha$ and
$\|\Qt{T}-(\lambda-B)^{-1}\|\le(C_1+\frac12)T^\alpha$, with
$C_1=4M+K^2\bigl(\frac{c_\alpha}8+\frac1{4\delta_1}\bigr)$. In particular the levels
of~\eqref{eq:mixed} are absolutely summable in operator norm.
\item[(ii)] Under \textup{(R$_\alpha$)}, for $x\in D(B)$,
\[
\|\Qtr{T}x-(\lambda-B)^{-1}x\|\le2K_R\,\zeta(1+\alpha)\Bigl(\frac T{2\pi}\Bigr)^{1+\alpha}
\|(\lambda-B)x\| .
\]
\item[(iii)] Under \textup{(D$_\alpha$)}, for $x\in X_\theta$, $0<\theta<1$,
$\|\Qtr{T}x-(\lambda-B)^{-1}x\|\le C\,T^{\alpha+\theta}\|x\|_{X_\theta}$.
\end{enumerate}
By Remark~\ref{rem:oneline}, \textup{(i)--(iii)} with $\alpha=1$ hold for every
analytic semigroup, and \textup{(ii)} then reads $\frac{K_R}{12}T^2\|(\lambda-B)x\|$.
\end{theorem}

\begin{proof}
(i) Let $x\in X$ and $g(t)=h(t)x$ with $h$ as in Lemma~\ref{lem:DtoR}. The trapezoid
error is $\sum_{m\ge0}e_m$ with
$e_m:=\frac T2\bigl(g(mT)+g((m+1)T)\bigr)-\int_{mT}^{(m+1)T}g$. Each panel satisfies
the crude bound $\|e_m\|\le2MT\|x\|$. For $m\ge1$, $g$ is $C^2$ on $[mT,(m+1)T]$, and the
Peano form of the trapezoid error,
$e_m=\frac12\int_{mT}^{(m+1)T}(t-mT)\bigl((m+1)T-t\bigr)g''(t)\,dt$ (valid for
$X$-valued functions by Hahn--Banach), gives
$\|e_m\|\le\frac{T^2}8\int_{mT}^{(m+1)T}\|g''\|$. Let $m_0:=\lceil T^{\alpha-1}\rceil$,
so that $T^\alpha\le m_0T\le T^\alpha+T\le2T^\alpha$. Using the crude bound for
$m<m_0$ and~\eqref{eq:h2int} with $\tau=m_0T\ge T^\alpha$ for $m\ge m_0$,
\[
\sum_{m\ge0}\|e_m\|\le2MT\,m_0\|x\|+\frac{K^2T^2}{8}\Bigl(\frac2{\delta_1}
+c_\alpha T^{\alpha(1-2/\alpha)}\Bigr)\|x\|
\le\Bigl(4M+\frac{c_\alpha K^2}8+\frac{K^2}{4\delta_1}\Bigr)T^\alpha\|x\| .
\]
The bound for $\Qt{T}=\Qtr{T}+\frac T2I$ follows from $T\le T^\alpha$. The level $k$
of~\eqref{eq:mixed} equals $\Qt{T_{k-1}}-\Qt{T_k}$, whose norm is at most
$(C_1+\frac12)(T_{k-1}^\alpha+T_k^\alpha)$ once $T_{k-1}\le1$, and
$\sum_kT_k^\alpha\le\sum_k(2^{-k}T_0)^\alpha<\infty$.

(ii) Insert \textup{(R$_\alpha$)} into~\eqref{eq:aliasbound}, noting $\nu_n\ge2\pi\ge1$:
$\sum_{n\ge1}2K_R\nu_n^{-1-\alpha}=2K_R\zeta(1+\alpha)(T/2\pi)^{1+\alpha}$. For
$\alpha=1$, $2\zeta(2)/(4\pi^2)=\frac1{12}$.

(iii) By Lemma~\ref{lem:DtoR}, \textup{(R$_\alpha$)} holds. For $y\in D(B)$, (i) and
(ii) give
$\|(\Qtr{T}-(\lambda-B)^{-1})x\|\le C_1T^\alpha\|x-y\|+C'T^{1+\alpha}\|y\|_{D(B)}
\le C\,T^\alpha\bigl(\|x-y\|+T\|y\|_{D(B)}\bigr)$; the infimum over $y$ is
$C\,T^\alpha K(T,x)\le C\,T^{\alpha+\theta}\|x\|_{X_\theta}$.
\end{proof}

Theorem~\ref{thm:analytic} is the analogue, for the sampled resolvent, of the
rate results for approximation formulas and rational approximations of semigroups
discussed in the introduction~\cite{GomilkoTomilov,GKT,BGT,Wakaiki24}; the proofs use
the same two ingredients, a time-domain estimate near $t=0$ and interpolation. All
three rates are sharp, already for normal generators.

\begin{proposition}[Sharpness]\label{prop:analytic-sharp}
Let $\lambda>0$.
\begin{enumerate}
\item[(a)] Let $A\ge0$ be self-adjoint with $\sigma(A)$ unbounded and $B=-A$. Then
\textup{(R$_1$)} holds with $K_R=1$, and for every $T>0$
\[
\|\Qtr{T}-(\lambda-B)^{-1}\|=\frac T2,\qquad\|\Qt{T}-(\lambda-B)^{-1}\|=T .
\]
If $0\in\sigma(A)$, then
$T^{-2}\|(\Qtr{T}-(\lambda-B)^{-1})(\lambda-B)^{-1}\|\to\frac1{12}$ as $T\to0^+$, so the
constant in Theorem~\textup{\ref{thm:analytic}(ii)} is attained.
\item[(b)] Let $0<\alpha\le1$, let $A$ have simple eigenvalues $a_j=2\pi2^j$,
$j\ge0$, with orthonormal eigenvectors $e_j$, and let $B_\alpha:=-iA-A^\alpha$. Then
$B_\alpha$ generates a contraction semigroup satisfying \textup{(D$_\alpha$)}, and
there are $c>0$ and $N_0$ such that for $T=T_N=2^{-N}$, $N\ge N_0$,
\[
\|\Qtr{T}-(\lambda-B_\alpha)^{-1}\|\ge cT^\alpha,\qquad
\bigl\|(\Qtr{T}-(\lambda-B_\alpha)^{-1})(\lambda-B_\alpha)^{-1}\bigr\|\ge cT^{1+\alpha},
\]
and, for $0<\theta<1$, the vector $x=\sum_j2^{-j\theta}e_j$ lies in $X_\theta$ and
satisfies $\|(\Qtr{T}-(\lambda-B_\alpha)^{-1})x\|\ge cT^{\alpha+\theta}$.
\end{enumerate}
\end{proposition}

\begin{proof}
(a) $\|(\lambda-i\nu+A)^{-1}\|=\sup_{a\in\sigma(A)}|\lambda+a-i\nu|^{-1}\le|\nu|^{-1}$.
With $u=T(\lambda+a)/2$ and the Langevin function $\mathcal L(u)=\coth u-1/u$, which
increases from $0$ to $1$ on $(0,\infty)$ and satisfies $\mathcal L(u)/u\to\frac13$ as
$u\to0^+$, the symbols of $\Qtr{T}-(\lambda-B)^{-1}$ and $\Qt{T}-(\lambda-B)^{-1}$ are
$\frac T2\mathcal L(u)$ and $\frac T2(1+\mathcal L(u))$; their suprema over the
unbounded set $\sigma(A)$ are $\frac T2$ and $T$. The symbol of
$(\Qtr{T}-(\lambda-B)^{-1})(\lambda-B)^{-1}$ is $\frac{T^2}4\mathcal L(u)/u$; at $a=0$
this tends to $T^2/12$, which with Theorem~\ref{thm:analytic}(ii) for $K_R=1$ gives the
limit.

(b) $B_\alpha$ is normal with eigenvalues $b_j=-ia_j-a_j^\alpha$, so $e^{tB_\alpha}$
is a contraction semigroup, and since $a_j^\alpha\le a_j$,
$\|(\lambda-B_\alpha)e^{tB_\alpha}\|\le\sup_{a\ge0}(\lambda+2a)e^{-ta^\alpha}
=\lambda+2(\alpha et)^{-1/\alpha}$, which is \textup{(D$_\alpha$)} with $\omega_1=0$.
The symbol of $\Qtr{T}-(\lambda-B_\alpha)^{-1}$ at $b_j$ is
$\varphi_T(p)=\frac T2\coth\frac{pT}2-\frac1p$ with $p=\lambda+a_j^\alpha+ia_j$. For
$j=N$ we have $a_NT=2\pi$, so $\coth\frac{pT}2=\coth u$ with
$u=\frac T2(\lambda+a_N^\alpha)>0$, and $\coth u\ge\max(1,1/u)$. Since $|p|\ge a_N$,
\[
|\varphi_T(p)|\ge\max\Bigl(\frac T2,\frac1{\lambda+(2\pi/T)^\alpha}\Bigr)-\frac{T}{2\pi}
\ge c\,T^\alpha\qquad(N\ge N_0),
\]
using the first entry if $\alpha=1$ and the second if $\alpha<1$. This gives the
first bound, and the second follows from $|p|\le(\lambda+2)a_N$. For the last claim,
the component $j=N$ contributes $|\varphi_T(p)|\,2^{-N\theta}\ge cT^{\alpha+\theta}$.
Finally, $D(B_\alpha)=D(A)$ with equivalent norms, and with $y$ the projection of $x$
onto $\operatorname{span}\{e_j:a_j\le1/t\}$ both $\|x-y\|$ and $t\|y\|_{D(A)}$ are
$O(t^\theta)$ for $0<t\le1$, so $K(t,x)=O(t^\theta)$ and $x\in X_\theta$.
\end{proof}

\begin{remark}[A scale of rates]\label{rem:dichotomy}
Theorems~\ref{thm:sa} and~\ref{thm:analytic} and Proposition~\ref{prop:riemann}
combine into one scale. With $X_0:=X$ and $X_1:=D(B)$, the error of $\Qtr{T}$ on
$X_s$ is of order $T^{\alpha+s}$, where $\alpha$ is the order of resolvent decay along
the alias line:
\begin{center}
\small
\begin{tabular}{lcccc}
\toprule
 & resolvent decay & operator norm & $x\in X_\theta$ & $x\in D(B)$\\
\midrule
unitary group, $\sigma(A)=\R$ & none ($\alpha=0$) & no convergence & $T^{\theta}$ & $T$\\
Crandall--Pazy, $0<\alpha<1$ & $|\nu|^{-\alpha}$ & $T^\alpha$ & $T^{\alpha+\theta}$ & $T^{1+\alpha}$\\
analytic semigroup & $|\nu|^{-1}$ & $T$ & $T^{1+\theta}$ & $T^2$\\
\bottomrule
\end{tabular}
\end{center}
All entries are sharp as bounds over the unit ball of the respective space, and for
$X_\theta$ they are attained by single vectors (Theorem~\ref{thm:sa}(iii), whose vector
lies in $X_{s'}$ for $0<s'<1$ by the $K$-functional argument at the end of the proof of
Proposition~\ref{prop:analytic-sharp}, and Proposition~\ref{prop:analytic-sharp}(b)). Each unit of
resolvent decay along the alias line buys one power of $T$. For $\alpha>0$ the
levels are absolutely summable, refining Remark~\ref{rem:order}. The
Euler--Maclaurin terms of Theorem~\ref{thm:EM} do not exist for $x\in X_\theta$, and
Romberg extrapolation does not improve on these rates. In the example of
Proposition~\ref{prop:analytic-sharp}(b), the mode $a_N$ sits on an alias of every
step $2^iT_N$, since $2^iT_Na_N\in2\pi\Z$. For $\alpha<1$ the symbols of all the
$\Qtr{2^iT_N}$ at that mode are asymptotic to $1/(\lambda+a_N^\alpha)$, and the
extrapolation weights sum to $1$, so every Romberg column of
Corollary~\ref{cor:romberg} has the same leading symbol there; for $\alpha=1$ the
columns retain a symbol of exact order $T_N$. The lower bounds of
Proposition~\ref{prop:analytic-sharp}(b) therefore persist.
\end{remark}

\begin{remark}[Resolvent decay alone]\label{rem:Ronly}
Theorem~\ref{thm:analytic}(ii) uses only \textup{(R$_\alpha$)}, and for $\alpha=1$ so do
(i) and (iii), by Remark~\ref{rem:oneline}. For $0<\alpha<1$, whether
\textup{(R$_\alpha$)} implies \textup{(D$_\alpha$)} depends on the space. On a
Hilbert space it does: \textup{(R$_\alpha$)} gives
$\|BS_t\|=O(t^{-1/\alpha})$ as $t\to0^+$~\cite[Thm.~2.3]{Wakaiki24},
\cite[Thm.~5.1]{Wakaiki25}, extending a result of Taylor~\cite{Taylor}, so all of
Theorem~\ref{thm:analytic} holds under \textup{(R$_\alpha$)}. On a Banach space,
\textup{(R$_\alpha$)} gives $\|BS_t\|=O(t^{1-2/\alpha})$~\cite{CrandallPazy}, improved
to $t^{-1/\alpha}$ up to a logarithmic factor~\cite[Thm.~3.1]{Wakaiki25}; the proof of (i) then gives the rate
$T^\alpha$ up to a logarithmic factor.
\end{remark}

%% =====================================================================
\section{Spectral measures}\label{sec:spectral}

Let $A$ be self-adjoint on $\Hil$, $\psi\in\Hil$, and let $\mu=\mu_\psi$ be its
spectral measure, a finite positive Borel measure with $\mu(\R)=\|\psi\|^2$. For
$E\in\R$ and $\eps>0$ we apply the constructions of \S\ref{sec:setting} to the
unitary group $e^{-it(A-E)}$ (with $\omega=0$, $M=1$, $\delta=\eps$) with Laplace parameter $\eps$,
writing $\Qt{T}(E,\eps)$ and $\F{T}{\zeta}(E,\eps)$. Since
$(\eps+i(A-E))^{-1}=-iR_A(E+i\eps)$ with $R_A(z)=(A-z)^{-1}$,
Theorem~\ref{thm:mixed} gives
\begin{equation}\label{eq:shifted}
R_A(E+i\eps)=i\Bigl[\Qt{T_0}(E,\eps)
-\sum_{k\ge1}T_k\sum_{r=1}^{b_k-1}\F{T_k}{\xi_k^r}(E,\eps)\Bigr].
\end{equation}
For $f\in C_c(\R)$ define the level functionals
\begin{align*}
L_0^\eps(f)&:=\frac1\pi\int_\R f(E)\,\Im\ip{i\Qt{T_0}(E,\eps)\psi}{\psi}\,dE,\\
L_k^\eps(f)&:=\frac1\pi\int_\R f(E)\,\Im\Bigl\langle iT_k\sum_{r=1}^{b_k-1}
\F{T_k}{\xi_k^r}(E,\eps)\psi,\psi\Bigr\rangle\,dE\qquad(k\ge1).
\end{align*}
The integrands are bounded continuous functions of $E$, so these are well defined.
By Stone's formula~\cite{RS1}, $\frac1\pi\int f(E)\Im\ip{R_A(E+i\eps)\psi}{\psi}dE
\to\int f\,d\mu$ as $\eps\to0^+$; the question is what the individual levels
contribute.

\begin{lemma}[Poisson-kernel form]\label{lem:poisson}
For $T>0$, $\eta\in\R$, $\zeta=e^{i\eta}$, $\eps>0$ and $E\in\R$,
\[
\Im\ip{iT\F{T}{\zeta}(E,\eps)\psi}{\psi}
=\frac T2\|\psi\|^2+\pi T\int_\R P_\varrho\bigl(T(a-E)-\eta\bigr)\,d\mu(a),
\qquad \varrho:=e^{-\eps T},
\]
where $P_\varrho(\phi)=\frac1{2\pi}\frac{1-\varrho^2}{1-2\varrho\cos\phi+\varrho^2}$ is
the Poisson kernel of the unit disc.
\end{lemma}

\begin{proof}
By the spectral theorem $\F{T}{\zeta}(E,\eps)$ acts as multiplication by $1/(1-w)$
with $w=\varrho e^{i(\eta-T(a-E))}$. Since $\Im(iz)=\Re z$ and
$\Re\frac1{1-w}=\frac12+\frac12\frac{1-|w|^2}{|1-w|^2}
=\frac12+\pi P_\varrho(\eta-T(a-E))$, and $P_\varrho$ is even, the claim follows.
\end{proof}

\begin{theorem}[Level functionals in closed form]\label{thm:spectral}
Let $A$ be self-adjoint, $\psi\in\Hil$, $f\in C_c(\R)$, and let $T_0$, $(b_k)$,
$\nu_0$, $\Lambda_k$ be as in \S\ref{sec:aliasing}. Then the limits
$L_k(f):=\lim_{\eps\to0^+}L_k^\eps(f)$ exist for all $k\ge0$, and
\begin{align}
L_0(f)&=\frac{T_0}{2\pi}\|\psi\|^2\!\int f+\sum_{\ell\in\Z}\int_\R f(a-\nu_0\ell)\,d\mu(a),
\label{eq:L0}\\
L_k(f)&=\frac{T_{k-1}-T_k}{2\pi}\|\psi\|^2\!\int f
+\sum_{\ell\in\Lambda_k}\int_\R f(a-\nu_0\ell)\,d\mu(a)\qquad(k\ge1),\label{eq:Lk}
\end{align}
where
$\sum_{\ell\in\Z}\int|f(a-\nu_0\ell)|\,d\mu(a)\le\|f\|_\infty\bigl(1+\diam(\supp f)/\nu_0\bigr)
\|\psi\|^2$. Moreover
\begin{equation}\label{eq:recovery}
\int_\R f\,d\mu=L_0(f)-\sum_{k=1}^\infty L_k(f),
\end{equation}
the series converging absolutely, and
$\lim_{\eps\to0^+}\sum_{k\ge1}L_k^\eps(f)=\sum_{k\ge1}L_k(f)$.
\end{theorem}

\begin{proof}
\emph{Step 1.} Fix $T>0$ and $\eta$. By Lemma~\ref{lem:poisson} and Fubini's
theorem (the kernel is nonnegative, $f$ is bounded with compact support, and $\mu$ is
finite),
\begin{gather*}
\frac1\pi\int f(E)\Im\ip{iT\F{T}{e^{i\eta}}(E,\eps)\psi}{\psi}dE
=\frac{T}{2\pi}\|\psi\|^2\!\int f+\int_\R\Phi_\eps(a)\,d\mu(a),\\
\Phi_\eps(a):=\int_\R f(E)\,TP_\varrho\bigl(T(a-E)-\eta\bigr)\,dE .
\end{gather*}
\emph{Step 2.} Substituting $E=a-(\eta+\phi)/T$ and using the $2\pi$-periodicity of
$P_\varrho$,
$\Phi_\eps(a)=\int_{-\pi}^{\pi}G(a-\phi/T)\,P_\varrho(\phi)\,d\phi$, where
$G(y):=\sum_{n\in\Z}f\bigl(y-(\eta+2\pi n)/T\bigr)$ is continuous and
$2\pi/T$-periodic, hence uniformly continuous, with
$\|G\|_\infty\le\|f\|_\infty(1+T\diam(\supp f)/2\pi)$. As $\varrho\to1^-$ the Poisson
integral converges uniformly, $\Phi_\eps\to G$, and $|\Phi_\eps|\le\|G\|_\infty$. By
dominated convergence,
$\int\Phi_\eps\,d\mu\to\int G\,d\mu=\sum_n\int f(a-(\eta+2\pi n)/T)\,d\mu(a)$, the
interchange being justified by Tonelli's theorem applied to $|f|$.

\emph{Step 3.} For $k=0$ take $T=T_0$, $\eta=0$; the frequencies are $\nu_0\ell$,
$\ell\in\Z$, which is~\eqref{eq:L0}. For $k\ge1$ sum Step~2 over
$\eta=2\pi r/b_k$, $1\le r<b_k$, with $T=T_k$; as in Corollary~\ref{cor:levels} the
frequencies are $\nu_0\Lambda_k$ and the constants add to
$(T_{k-1}-T_k)/2\pi$, which is~\eqref{eq:Lk}. The bound on
$\sum_\ell\int|f(a-\nu_0\ell)|\,d\mu$ holds because for each $a$ at most
$1+\diam(\supp f)/\nu_0$ indices $\ell$ have $a-\nu_0\ell\in\supp f$.

\emph{Step 4.} The $\Lambda_k$ are disjoint, so
$\sum_k|L_k(f)|\le\frac{T_0}{2\pi}\|\psi\|^2\int|f|
+\sum_{\ell\ne0}\int|f(a-\nu_0\ell)|\,d\mu<\infty$. Since they partition $\Z\setminus\{0\}$
and $\sum_{k\ge1}(T_{k-1}-T_k)=T_0$, subtracting~\eqref{eq:Lk} summed over $k$
from~\eqref{eq:L0} leaves the term $\ell=0$, which is~\eqref{eq:recovery}.

\emph{Step 5.} For fixed $\eps>0$, the partial sums of~\eqref{eq:shifted} are
$i\Qt{T_N}(E,\eps)$, bounded in norm by $1/\eps+T_0$ uniformly in $E$ and $N$
(Lemma~\ref{lem:bound}). By dominated convergence in $E$,
$\frac1\pi\int f\Im\ip{R_A(E+i\eps)\psi}{\psi}dE=L_0^\eps(f)-\sum_kL_k^\eps(f)$. Letting
$\eps\to0^+$ and using Stone's formula and Step~3 for $k=0$,
$\sum_kL_k^\eps(f)\to L_0(f)-\int f\,d\mu=\sum_kL_k(f)$.
\end{proof}

\begin{remark}[Scope]\label{rem:scope}
Theorem~\ref{thm:spectral} holds for every self-adjoint operator and every vector:
no assumption is made on the spectral type of $\mu_\psi$ (point, absolutely or
singular continuous, embedded eigenvalues), on its support, or on $f$ beyond
$f\in C_c(\R)$. The recovery~\eqref{eq:recovery} of $\int f\,d\mu$ is, of course, Stone's
formula, which needs no hypotheses either; the content of
Theorem~\ref{thm:spectral} is the accounting of what each level contributes.
In~\cite[Prop.~8, Cor.~4, Prop.~9]{CastilloV1} the corresponding statement was proved
for $A=-\Delta$ with $\widehat\psi$ supported in an annulus, and for $A=-\Delta+V$
under a limiting absorption principle, with $f\in C^1_c((0,\infty))$; for
$b_k\equiv2$ and $T_0=1$ the functionals $L_k$ coincide with the $J_k$ there. Those
hypotheses are unnecessary: the limiting absorption principle concerns the existence
of a continuous density for $\mu_\psi$ and plays no role in the recovery.
\end{remark}

\begin{remark}[What the levels do]\label{rem:nyquist}
The constants
$c_0(f):=\frac{T_0}{2\pi}\|\psi\|^2\int f$ and
$c_k(f):=\frac{T_{k-1}-T_k}{2\pi}\|\psi\|^2\int f$ are known, and they telescope. Writing
$\widehat L_k:=L_k-c_k$ for the alias parts,
$\int f\,d\mu=\widehat L_0(f)-\sum_{k\ge1}\widehat L_k(f)$.
\begin{enumerate}
\item[(a)] If $(\supp\mu-\supp f)\cap\nu_0(\Z\setminus\{0\})=\emptyset$, then
$\widehat L_k(f)=0$ for all $k\ge1$ and $\int f\,d\mu=\widehat L_0(f)$: a single sampling
rate above the Nyquist rate of the spectral window suffices, and the dyadic levels
are not needed.
\item[(b)] In general, an index $\ell\in\Lambda_k$ has $|\ell|\ge\beta_{k-1}$, so
\[
|\widehat L_k(f)|\le\|f\|_\infty\Bigl(1+\frac{\diam(\supp f)}{\nu_0}\Bigr)\,
\mu\bigl(\{a:\operatorname{dist}(a,\supp f)\ge\nu_0\beta_{k-1}-\diam(\supp f)\}\bigr).
\]
The levels are needed exactly when the spectral support of $\psi$, seen from the
window $\supp f$, extends beyond the Nyquist band, and the size of level $k$ is
controlled by the mass of $\mu_\psi$ at distance $\approx\nu_0\beta_{k-1}$.
\end{enumerate}
The uncorrected partial sums $L_0-\sum_{k\le K}L_k$ converge only like $T_K$, because
of the constants; the corrected sums $\widehat L_0-\sum_{k\le K}\widehat L_k$ are exact
as soon as the remaining aliases carry no mass (Table~\ref{tab:recon}).
\end{remark}

%% =====================================================================
\section{Numerical illustrations}\label{sec:numerics}

All computations use double precision; the scripts are available as ancillary files.
Unless stated otherwise $b_k\equiv2$, $T_0=1$ and $T_N=2^{-N}$.

\subsection{Rates for a self-adjoint generator}\label{sec:num-rates}
\begin{figure}[ht]
\centering
\includegraphics[width=\textwidth]{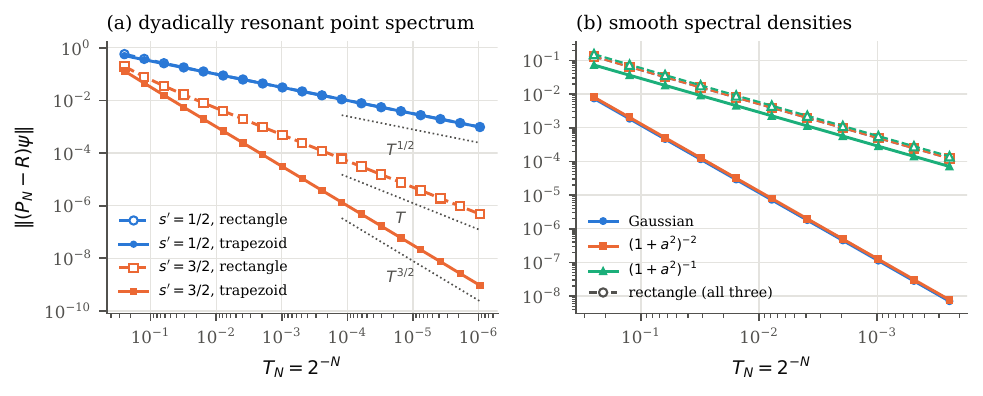}
\caption{Error of the partial sums of~\eqref{eq:CCC} (rectangle, $\Qt{T_N}$) and of
the trapezoid correction ($\Qtr{T_N}$) for $A$ = multiplication by $a$, $\lambda=1$.
(a)~Spectral mass on the aliases $2\pi2^j$, weights $\propto2^{-2js'}$: the rates
$T_N^{\min(s',1)}$ and $T_N^{\min(s',2)}$ of Theorem~\ref{thm:sa} are attained
(dotted guides). (b)~Smooth spectral densities: faster rates, cf.\
Remark~\ref{rem:smooth}; the three rectangle curves nearly coincide.}
\label{fig:rates}
\end{figure}

Let $A$ be multiplication by $a$, so that
$\|(\Qt{T}-(\lambda+iA)^{-1})\psi\|^2=\int|e_T(a)|^2d\mu_\psi(a)$ with $e_T$ as in the
proof of Theorem~\ref{thm:sa}, and $\lambda=1$. Figure~\ref{fig:rates}(a) uses the
dyadically resonant point spectrum of Theorem~\ref{thm:sa}(iii), with weights
$\mu_\psi(\{2\pi2^j\})\propto2^{-2js'}$; the fitted orders over the last four points are
$0.500$ for $s'=\frac12$ (both rules) and $1.000$ (rectangle), $1.500$ (trapezoid)
for $s'=\frac32$, as predicted. Figure~\ref{fig:rates}(b) uses absolutely continuous
measures with densities $e^{-a^2/2}$, $(1+a^2)^{-2}$ and $(1+a^2)^{-1}$, for which
$\psi\in D(|A|^s)$ for all $s$, for $s<\frac32$ and for $s<\frac12$ respectively. The
rectangle rule converges at order $1.000$ in all three cases; the trapezoid rule at
orders $2.000$, $2.000$ and $1.000$, in agreement with Remark~\ref{rem:smooth}.

\subsection{The square-lattice Green function}\label{sec:num-lattice}
For the adjacency operator $H$ of $\Z^2$, the return Green function
$G(z)=\ip{\delta_0}{(z-H)^{-1}\delta_0}$ has the closed form
$G(z)=\frac{2}{\pi z}K(4/z)$~\cite{Economou}, and
$\ip{\delta_0}{e^{-itH}\delta_0}=J_0(2t)^2$ (cf.~\cite[Example~3]{CastilloV1}). The
partial sum through level $N$ is
\[
G_N=-iT_N\sum_{m\ge0}e^{izmT_N}J_0(2mT_N)^2 .
\] At $z=1.3+0.2i$ we checked the closed
form against direct quadrature of $-i\int_0^\infty e^{izt}J_0(2t)^2dt$ to double
precision. Figure~\ref{fig:accel}(a) shows the error. The rectangle rule has error
$1.22\times10^{-4}=T_{12}/2$ at $N=12$, as Corollary~\ref{cor:error} predicts when the
aliases are negligible; and the trapezoid correction has error $6.5\times10^{-9}$. The
Romberg columns $P^{(1)}_N$ and $P^{(2)}_N$ shown in the figure reach rounding level
at $N=12$ and $N=9$, and the next column $P^{(3)}_N$ (not shown) at $N=7$; at $N=9$ and
$N=7$ this takes about $1.2\times10^5$ and $2.9\times10^4$ propagator samples
respectively. The uncorrected series would need $N\approx52$ for the same accuracy.

\subsection{A non-normal semigroup}\label{sec:nonnormal}
\begin{figure}[ht]
\centering
\includegraphics[width=\textwidth]{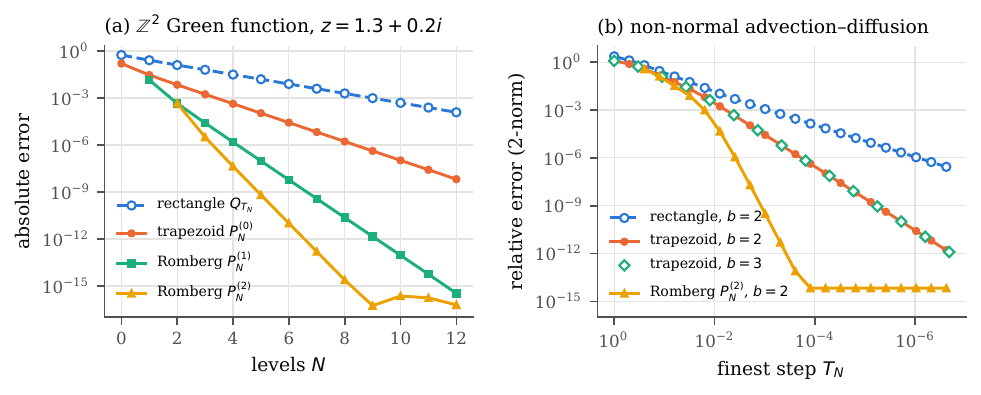}
\caption{Acceleration. (a)~Return Green function of $\Z^2$ at $z=1.3+0.2i$: rectangle
rule, trapezoid correction, and Romberg columns $P^{(1)}_N$, $P^{(2)}_N$
(Corollary~\ref{cor:romberg}).
(b)~Resolvent of a non-normal advection--diffusion matrix at $\lambda=1$: rectangle and
trapezoid rules for radix $b=2$ and $b=3$, and the second Romberg column for $b=2$.}
\label{fig:accel}
\end{figure}

Let $B=\nu\partial_x^2-c\,\partial_x$ on $(0,1)$ with Dirichlet conditions, discretized by
centred differences on $n=60$ interior points, with $\nu=0.01$ and $c=1$. The matrix is
tridiagonal Toeplitz, with real eigenvalues
$-2\nu/h^2+2\bigl(\nu^2/h^4-c^2/(4h^2)\bigr)^{1/2}\cos(k\pi h)$, $h=1/(n+1)$,
$1\le k\le n$, so its spectral abscissa is $-31.85$; it has $\|B\|_2=148.8$, and its
eigenvector matrix has condition number above $10^{20}$. It is therefore far from
normal, and Theorem~\ref{thm:sa} does not apply~\cite{TE}; a double-precision
eigensolver returns the abscissa $-29.04$. With $\lambda=1$,
Figure~\ref{fig:accel}(b) shows the relative error in the $2$-norm of the partial
sums~\eqref{eq:partial} computed level by level. While $T_N\|B\|\gtrsim1$ both rules
converge at first order in operator norm, as Theorem~\ref{thm:analytic}(i) predicts
for the analytic semigroup generated by the continuous operator (see
\S\ref{sec:num-analytic}); once $T_N\|B\|\lesssim1$ the matrix regime of
Theorem~\ref{thm:EM} takes over, the rectangle rule staying first order and the
trapezoid correction becoming second order, for $b=2$ and $b=3$ alike. As
Theorem~\ref{thm:EM} predicts, the error depends on $T_N$ only. The Romberg column $P^{(2)}_N$ reaches
$7\times10^{-15}$ at $T_N=2^{-13}$. Evaluating $\Qtr{T_N}$ directly rather than through
the levels gives the same errors down to $T_N=2^{-20}$, but at $T_N=2^{-22}$ the
direct error is $1.4\times10^{-10}$ against $1.6\times10^{-12}$ for the series
(Remark~\ref{rem:stability}).

\subsection{Analytic versus unitary}\label{sec:num-analytic}
\begin{figure}[ht]
\centering
\includegraphics[width=\textwidth]{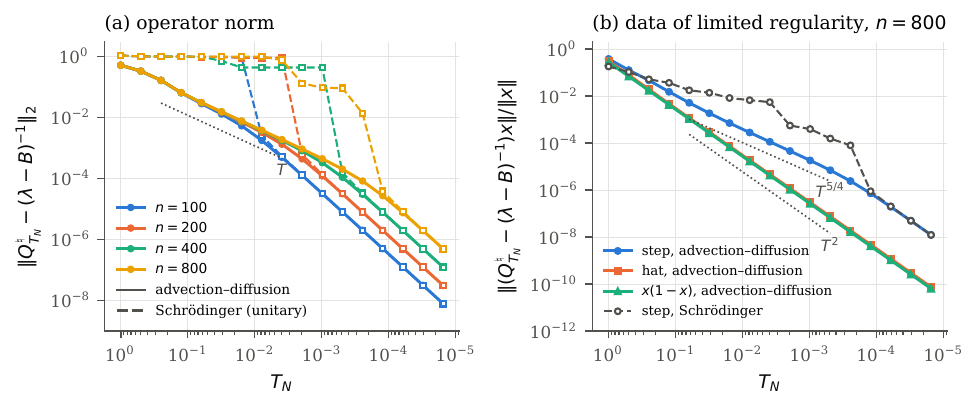}
\caption{Analytic versus unitary. (a)~Operator-norm error of $\Qtr{T_N}$ for the
advection--diffusion matrix $B_a$ (solid) and the Schr\"odinger matrix $B_u$ (dashed),
$n=100$--$800$: uniform first-order convergence for the analytic semigroup, an
$n$-dependent plateau for the unitary group. (b)~Errors for data of limited
regularity, $n=800$; dotted guides $T^{5/4}$ and $T^2$.}
\label{fig:analytic}
\end{figure}

To test Theorem~\ref{thm:analytic} and Remark~\ref{rem:dichotomy} independently of
the discretization, we compare $B_a=\nu D_2-cD_1$, the advection--diffusion matrix
of \S\ref{sec:nonnormal}, with $B_u=i\nu D_2$, which generates a unitary group (a
discrete Schr\"odinger equation), for $n=100,200,400,800$, $\nu=0.01$, $c=1$,
$\lambda=1$. Both matrices have $\|B\|_2$ growing like $n^2$, from $408$ to
$2.6\times10^4$. Figure~\ref{fig:analytic}(a) shows the operator-norm error of the
trapezoid-corrected partial sums. For $B_a$ the curves coincide for all $n$ and decay
at first order (local orders $1.00$--$1.08$ for $n=800$ and $3\le N\le10$) until
$T_N\|B_a\|\lesssim1$, after which the finite-dimensional second-order regime begins:
the bound of Theorem~\ref{thm:analytic}(i) is uniform in the discretization. For
$B_u$ the error stays of order one until $2\pi/T_N$ is comparable to $\|B_u\|$, i.e.\
until the first alias leaves the spectrum, so the onset of convergence recedes as $n$
grows and there is no convergence in operator norm in the continuum limit.

Figure~\ref{fig:analytic}(b) uses $n=800$ and three vectors: the indicator of
$(\frac14,\frac34)$, the hat function $\min(x,1-x)$, and $x(1-x)$. Their sine
coefficients decay like $k^{-1}$, $k^{-2}$ and $k^{-3}$; since $D(B_a)$ coincides with
the domain of the Dirichlet Laplacian, a $K$-functional computation in the sine basis
places them in $X_{1/4}$, $X_{3/4}$ and $D(B_a)$. For $x(1-x)$ the order is $2.00$, as
in Theorem~\ref{thm:analytic}(ii). For the indicator the local order falls to $1.30$
at $N\approx9$, approaching the predicted $1+\theta=1.25$, before the matrix regime
restores order $2$. For the hat function the observed order is $2.0$ throughout: the
$T^{7/4}$ contribution of Theorem~\ref{thm:analytic}(iii) comes from high modes and has
a small constant, and in the accessible range it is dominated by the $T^2$
contribution of the low modes. The indicator under the unitary group $B_u$ converges
slowly and irregularly (overall order $\approx0.65$ for $N\le8$), as the aliases pass
through the spectrum.

For the diagonal example of Proposition~\ref{prop:analytic-sharp}(a) the predicted
quantities are reproduced exactly: $\|\Qtr{T}-(\lambda-B)^{-1}\|/(T/2)=1.000000$ and
$T^{-2}\|(\Qtr{T}-(\lambda-B)^{-1})(\lambda-B)^{-1}\|\cdot12=0.999998$ at $T=10^{-2}$.
For $B=-A$ with eigenvalues $2^j$ and $x=\sum_j2^{-j\theta}e_j\in X_\theta$ the fitted
orders are $1.250$, $1.500$, $1.750$ for $\theta=\frac14,\frac12,\frac34$.
We also checked~\eqref{eq:aliasline} on a random non-normal $8\times8$ matrix: the
discrepancy after $2000$ alias pairs equals the predicted tail.

\subsection{Intermediate decay: a structurally damped string}\label{sec:num-damped}
To test the intermediate range $0<\alpha<1$ of Theorem~\ref{thm:analytic}, i.e.\ the
Crandall--Pazy class with parameter $\alpha<1$, we take the damped string $\ddot u+\rho A^\beta\dot u+Au=0$, with $A=-\partial_x^2$ on $(0,1)$ with
Dirichlet conditions, $\rho=1$ and $\beta\in\{0,\frac18,\frac14,\frac38,\frac12\}$, in
the energy space $D(A^{1/2})\times L^2$. In modal energy coordinates the generator is
block diagonal, with blocks
\[
B_k=\begin{pmatrix}0&s_k\\-s_k&-d_k\end{pmatrix},\qquad s_k=\pi k,\quad d_k=\rho s_k^{2\beta}.
\]
Each block has eigenvalues $-d_k/2\pm iw_k$, $w_k=(s_k^2-d_k^2/4)^{1/2}$, and an
eigenvector matrix with condition number at most $\sqrt3$. Hence
\textup{(R$_{2\beta}$)} and \textup{(D$_{2\beta}$)} follow from the scalar estimates
$\sup_k|\lambda-i\nu+d_k/2\mp iw_k|^{-1}\lesssim|\nu|^{-2\beta}$ and
$\sup_k(\lambda+s_k)e^{-td_k/2}\lesssim t^{-1/(2\beta)}$, and the theory predicts
$\alpha=2\beta$. For $\beta=0$ the spectrum lies on the line $\Re z=-\rho/2$ and
$\alpha=0$; for $\beta=\frac12$ the semigroup is analytic.

The error operator $\Qtr{T}-(\lambda-B)^{-1}$ is block diagonal with blocks
$\varphi_T(\lambda-B_k)$, $\varphi_T(p)=\frac T2\coth\frac{pT}2-\frac1p$, which we
evaluate by two-point interpolation on the spectrum of each block. Its operator norm
is the supremum over $k$. We keep $2^{21}$ modes, whose frequencies reach the
$16$th alias at the finest step $T=2^{-16}$, and take $\lambda=1$. Against a
$40$-digit reference this evaluation has relative error below $10^{-11}$; computing
$\Qt{T}$ directly from $e^{TB_k}$ loses about eight digits to cancellation. By
\eqref{eq:partial} these are exactly the partial sums of the dyadic series; the
level-by-level evaluation was tested in \S\ref{sec:nonnormal}.

\begin{figure}[ht]
\centering
\includegraphics[width=\textwidth]{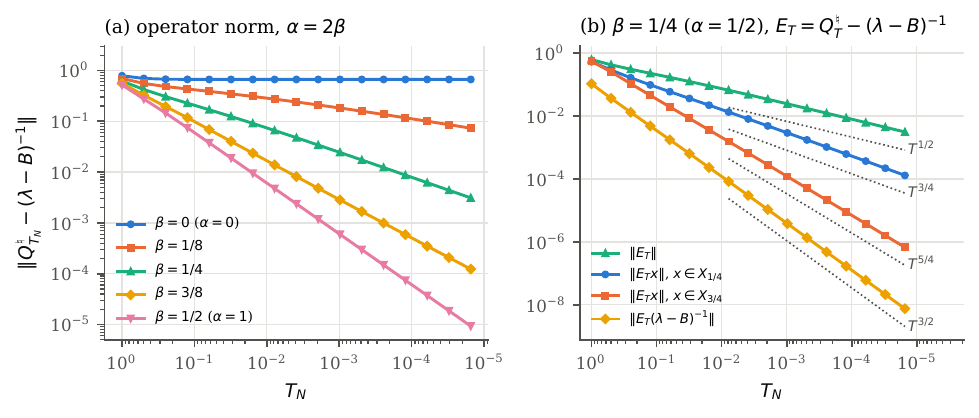}
\caption{Structurally damped string. (a)~Operator-norm error of $\Qtr{T_N}$ for
damping exponents $\beta=0,\dots,\frac12$; the order is $\alpha=2\beta$ (Table~\ref{tab:damped}).
(b)~$\beta=\frac14$: operator norm, norm on $D(B)$, and resonant vectors in $X_{1/4}$
and $X_{3/4}$, with the sharp rates of Theorem~\ref{thm:analytic} as dotted guides.}
\label{fig:damped}
\end{figure}

Table~\ref{tab:damped} and Figure~\ref{fig:damped} show the results. In every case
the operator-norm order agrees with the measured exponent of resolvent decay to
within $0.01$, and the order on $D(B)$, measured through
$\|(\Qtr{T}-(\lambda-B)^{-1})(\lambda-B)^{-1}\|$, exceeds it by one, as
Theorem~\ref{thm:analytic} predicts. For $\beta=0$ the error stays at
$1/(\lambda+\rho/2)=\frac23$. For $\beta=\frac18$ both exponents are $0.215$ instead
of $\frac14$, a pre-asymptotic effect: at the first alias the error is
$\approx1/(\lambda+\rho s^{1/4}/2)$ with $\rho s^{1/4}/2\approx12.6$ at $N=16$, whose
local exponent is $\frac14\cdot\frac{12.6}{13.6}\approx0.23$, matching the observed
local order $0.230$. The time-domain bound \textup{(D$_{2\beta}$)} is confirmed by the
slopes of $\|(\lambda-B)S_t\|$: $-1.95$, $-1.30$ and $-0.99$ for
$\beta=\frac14,\frac38,\frac12$, against $-2$, $-\frac43$ and $-1$. For $\beta=\frac14$,
the vectors supported on the modes $k=2^{j+1}$, whose frequencies sit at the aliases
$2\pi2^j$, with weights $2^{-j\theta}$, lie in $X_\theta$ and converge with orders
$0.745$ and $1.243$ for $\theta=\frac14$ and $\frac34$, against the sharp
$\alpha+\theta=0.75$ and $1.25$ of Theorem~\ref{thm:analytic}(iii) and
Proposition~\ref{prop:analytic-sharp}(b).

\begin{table}[ht]
\centering
\small
\begin{tabular}{cccccc}
\toprule
$\beta$ & $\alpha=2\beta$ & resolvent decay & operator norm & on $D(B)$ & predicted\\
\midrule
$0$ & $0$ & $0.000$ & $0.000$ & $1.000$ & $\alpha$, $1+\alpha$\\
$1/8$ & $0.25$ & $0.215$ & $0.215$ & $1.215$ & \\
$1/4$ & $0.5$ & $0.492$ & $0.492$ & $1.494$ & \\
$3/8$ & $0.75$ & $0.749$ & $0.755$ & $1.764$ & \\
$1/2$ & $1$ & $1.000$ & $1.000$ & $2.000$ & \\
\bottomrule
\end{tabular}
\caption{Damped string: fitted exponents of $\|(\lambda-i\nu-B)^{-1}\|$ over
$2^{10}\le\nu\le2^{20}$ (sign reversed), and fitted orders of
$\|\Qtr{T}-(\lambda-B)^{-1}\|$ and $\|(\Qtr{T}-(\lambda-B)^{-1})(\lambda-B)^{-1}\|$ over
$2^{-16}\le T\le2^{-8}$.}
\label{tab:damped}
\end{table}

\subsection{Spectral measures}\label{sec:num-spectral}
Let $\mu=0.3\,\delta_{\{1\}}+0.7\,\mathcal N(0.5,6^2)$, the sum of an eigenvalue
($\delta_{\{1\}}$ is the unit point mass at $1$) and an absolutely continuous part
with unbounded support, and
$f(E)=\exp(-1/(1-u^2))$ for $|u|<1$, $u=(E-1)/2$, so $\supp f=[-1,3]$. With $T_0=1$ the
Nyquist band has width $\nu_0=2\pi$, and the Gaussian tail produces aliases at
levels $0$ through $3$. Table~\ref{tab:levels} compares $L_k^\eps(f)$, computed by
quadrature directly from the resolvent symbols, with the closed forms of
Theorem~\ref{thm:spectral}; the differences are $O(\eps)$. Table~\ref{tab:recon} shows
the reconstruction~\eqref{eq:recovery}: the uncorrected sums converge like $2^{-K}$,
while the corrected sums of Remark~\ref{rem:nyquist} are exact to rounding after $K=3$
levels, because the next aliases, at distance $16\pi$, carry Gaussian mass of order
$10^{-16}$.

\begin{table}[ht]
\centering
\small
\begin{tabular}{ccccc}
\toprule
$k$ & $L_k^\eps(f)$, $\eps=10^{-3}$ & $L_k^\eps(f)$, $\eps=10^{-4}$ & $L_k(f)$, closed form & alias part $\widehat L_k(f)$\\
\midrule
0 & 0.35055338 & 0.35061405 & 0.35062079 & $2.09\times10^{-1}$\\
1 & 0.11915152 & 0.11914611 & 0.11914550 & $4.85\times10^{-2}$\\
2 & 0.04493838 & 0.04493388 & 0.04493338 & $9.60\times10^{-3}$\\
3 & 0.01768279 & 0.01768164 & 0.01768151 & $1.56\times10^{-5}$\\
4 & 0.00883326 & 0.00883300 & 0.00883298 & $<10^{-15}$\\
5 & 0.00441656 & 0.00441650 & 0.00441649 & $0$\\
\bottomrule
\end{tabular}
\caption{Level functionals for $\mu=0.3\,\delta_{\{1\}}+0.7\,\mathcal N(0.5,6^2)$: direct
quadrature at small $\eps$ versus the closed forms~\eqref{eq:L0}--\eqref{eq:Lk}.}
\label{tab:levels}
\end{table}

\begin{table}[ht]
\centering
\small
\begin{tabular}{ccc}
\toprule
$K$ & error of $L_0-\sum_{k\le K}L_k$ & error of $\widehat L_0-\sum_{k\le K}\widehat L_k$\\
\midrule
0 & $1.99\times10^{-1}$ & $5.81\times10^{-2}$\\
1 & $8.03\times10^{-2}$ & $9.62\times10^{-3}$\\
2 & $3.53\times10^{-2}$ & $1.56\times10^{-5}$\\
3 & $1.77\times10^{-2}$ & $1.1\times10^{-16}$\\
5 & $4.42\times10^{-3}$ & $0$\\
10 & $1.38\times10^{-4}$ & $0$\\
20 & $1.35\times10^{-7}$ & $0$\\
\bottomrule
\end{tabular}
\caption{Reconstruction of $\int f\,d\mu=0.15119443948767$ from the closed-form
levels, without and with removal of the known constants $c_k(f)$.}
\label{tab:recon}
\end{table}

%% =====================================================================
\section{Green functions from sampled autocorrelations}\label{sec:practical}

In many applications the propagator of a self-adjoint Hamiltonian $H$ is accessible
only through an autocorrelation record
\[
C(t)=\ip{e^{-itH}\psi}{\psi}=\int_\R e^{-ita}\,d\mu_\psi(a),
\]
sampled on a uniform clock $t_m=m\Delta t$: real-time electronic-structure and
wave-packet propagation, molecular dynamics, experiment, and Hamiltonian simulation.
The quantity sought is the Green function or the broadened spectral density,
\[
G(z)=\ip{(z-H)^{-1}\psi}{\psi}=\int_\R\frac{d\mu_\psi(a)}{z-a}
=-i\int_0^\infty e^{izt}C(t)\,dt,\qquad z=E+i\eta,\ \eta>0,
\]
and $-\pi^{-1}\Im G(E+i\eta)$. This section restates the results of
Sections~\ref{sec:aliasing}--\ref{sec:spectral} in that language. We set $\hbar=1$;
otherwise $2\pi/\Delta t$ reads $2\pi\hbar/\Delta t$ throughout.

\emph{Dictionary.} Put $B=-iH$, $\lambda=-iz=\eta-iE$ and $T=\Delta t$, so that
\textup{(S)} holds with $M=1$, $\omega=0$, $\delta=\Re\lambda=\eta$, and
$\lambda-B=i(H-z)$. Then $G(z)=-i\ip{(\lambda-B)^{-1}\psi}{\psi}$, and the rectangle
and trapezoid sums of the record are
\begin{align*}
G_{\Delta t}(z)&:=-i\Delta t\sum_{m\ge0}e^{izm\Delta t}C(m\Delta t)
=-i\ip{\Qt{\Delta t}\psi}{\psi},\\
G^\natural_{\Delta t}(z)&:=G_{\Delta t}(z)+\frac{i\Delta t}{2}\,C(0)
=-i\ip{\Qtr{\Delta t}\psi}{\psi}.
\end{align*}
The alias points of Theorem~\ref{thm:alias} are $\lambda-2\pi in/\Delta t
=-i(z+2\pi n/\Delta t)$, so for $\psi\in D(H)$
\begin{equation}\label{eq:periodized}
G^\natural_{\Delta t}(z)=\lim_{L\to\infty}\sum_{|n|\le L}
G\Bigl(z+\frac{2\pi n}{\Delta t}\Bigr):
\end{equation}
the trapezoid sum is the Green function periodized in energy with period
$2\pi/\Delta t$, and the rectangle sum adds the bias $-\frac{i\Delta t}{2}C(0)$, which is
the leading error for small $\Delta t$ and should always be removed.

\begin{proposition}[Explicit error bound]\label{prop:practical}
Let $\psi\in\Hil$, $z=E+i\eta$ with $\eta>0$, $\Delta t>0$, and let
$W:=\{a\in\R:\ |a-z|\le\pi/\Delta t\}$. Then
\[
\bigl|G^\natural_{\Delta t}(z)-G(z)\bigr|
\le\frac{\Delta t^2}{9}\int_W|a-z|\,d\mu_\psi(a)
+\Bigl(\frac2\eta+\frac{3\Delta t}{2}\Bigr)\mu_\psi(\R\setminus W).
\]
\end{proposition}

\begin{proof}
By the spectral theorem, $G^\natural_{\Delta t}(z)-G(z)=-i\int e(a)\,d\mu_\psi(a)$ with
$e(a)=\Delta t/(1-e^{-\Delta t\,p})-\Delta t/2-1/p$ and $p=\eta+i(a-E)$, so that
$|p|=|a-z|$. Since $\Delta t/(1-e^{-\Delta t\,p})-\Delta t/2=\frac{\Delta t}2
\coth\frac{\Delta t\,p}2$, the partial-fraction expansion of $\coth$ gives
$e(a)=\sum_{n\ge1}2p/(p^2+\nu_n^2)$ with $\nu_n=2\pi n/\Delta t$, the scalar case of
Theorem~\ref{thm:aliasline}. For $a\in W$ we have $|p|\le\nu_1/2$, hence
$|p^2+\nu_n^2|\ge\nu_n^2-|p|^2\ge\frac34\nu_n^2$ and
$|e(a)|\le\frac83|p|\sum_{n\ge1}\nu_n^{-2}=\frac83|p|\frac{\Delta t^2}{24}
=\frac{\Delta t^2}{9}|p|$. For $a\notin W$, $|1-e^{-\Delta t\,p}|\ge1-e^{-\eta\Delta t}$
and $x/(1-e^{-x})\le1+x$ for $x>0$ give
$|e(a)|\le(\eta^{-1}+\Delta t)+\Delta t/2+\eta^{-1}$.
\end{proof}

The following consequences are the practical content of the paper for this setting.

\begin{enumerate}
\item[(a)] \emph{Step size.} If $\mu_\psi$ is supported in a bounded interval $J$ (a
finite basis, a lattice, or $\psi$ filtered to an energy band), the choice
\[
\Delta t\le\frac{\pi}{\max_{a\in J}|a-z|}
\]
puts the whole spectrum of $\psi$ in the window $W$ around the target energy. The
second term of Proposition~\ref{prop:practical} then vanishes, and the error is at
most $\frac{\Delta t^2}{9}\max_J|a-z|\,\|\psi\|^2$, with no dependence on $\eta$. When
a range of energies $E$ is wanted, the condition must hold for each of them. If the
step is coarser, spectral weight outside $W$ is aliased into the window
by~\eqref{eq:periodized}, and the error can be of order
$\mu_\psi(\R\setminus W)/\eta$, independently of $\Delta t$, until that weight has
left the aliased bands. For unbounded spectrum with $\int|a|^{2s}d\mu_\psi<\infty$,
Chebyshev's inequality bounds the off-window mass by
$(\Delta t/\pi)^{2s}\int|a-z|^{2s}d\mu_\psi$, and the error is
$O(\Delta t^{\min(2,2s)})$. This exponent cannot be improved: for the resonant
spectral measure of Theorem~\ref{thm:sa}(iii), with $s'<1$, the error is of exact order
$\Delta t^{2s'}$ while $\int|a|^{2s}d\mu_\psi<\infty$ for every $s<s'$.

\item[(b)] \emph{Moment corrections.} With $\lambda-B=i(H-z)$,
Theorem~\ref{thm:EM} reads, for $\psi\in D(H^{2q})$,
\[
G^\natural_{\Delta t}(z)=G(z)+\sum_{j=1}^{q-1}\frac{|B_{2j}|}{(2j)!}\,\Delta t^{2j}
\ip{(H-z)^{2j-1}\psi}{\psi}+r_q,\qquad
|r_q|\le\frac{2|B_{2q}|}{(2q)!\,\eta}\,\Delta t^{2q}\,\|(H-z)^{2q}\psi\|\,\|\psi\| .
\]
The first two terms are $\frac{\Delta t^2}{12}\bigl(\ip{H\psi}{\psi}-z\|\psi\|^2\bigr)$
and $\frac{\Delta t^4}{720}\ip{(H-z)^3\psi}{\psi}$. The moments
$\ip{H^k\psi}{\psi}=i^kC^{(k)}(0)$ are often known exactly (sum rules, or a few
applications of $H$ to $\psi$), and then these terms can be subtracted as
end corrections of Gregory type, raising the order to $\Delta t^4$ or $\Delta t^6$
at no cost in samples.

\item[(c)] \emph{Romberg from one record.} The sum $G^\natural_{2\Delta t}$ uses only the
samples $C(2m\Delta t)$, so a single record at step $\Delta t$ yields the whole
Romberg table of Corollary~\ref{cor:romberg} with $b=2$. This is the moment-free
alternative to~(b), under the same smoothness hypothesis.

\item[(d)] \emph{Record length.} Since $|C(t)|\le\|\psi\|^2$, truncating the sum at
$t_{\max}=m_{\max}\Delta t$ costs at most
$\Delta t\,e^{-\eta t_{\max}}\|\psi\|^2/(1-e^{-\eta\Delta t})
\le(\eta^{-1}+\Delta t)\,e^{-\eta t_{\max}}\|\psi\|^2$, so $t_{\max}$ of a few tens of
$\eta^{-1}$ suffices for double precision. The broadening $\eta$ sets the energy
resolution, and halving it doubles the record.

\item[(e)] \emph{Integrated spectral quantities.} For $\int f\,d\mu_\psi$ with $f$ a
continuous spectral window, obtained from Stone's formula as $\eta\to0^+$,
Theorem~\ref{thm:spectral} and Remark~\ref{rem:nyquist} give the aliased part of each
level in closed form. With $T_0=\Delta t$: if no difference between a point of
$\supp f$ and a point of $\supp\mu_\psi$ is a nonzero multiple of $2\pi/\Delta t$,
the single step $\Delta t$ already gives $\int f\,d\mu_\psi$ up to the known constant
$c_0(f)$, and no refinement is needed; otherwise the levels remove the aliased
weight at the rate described in Remark~\ref{rem:nyquist}(b).
\end{enumerate}

For the lattice Green function of \S\ref{sec:num-lattice}, $\supp\mu=[-4,4]$ and
$z=1.3+0.2i$, so~(a) requires $\Delta t\le\pi/|{-4}-z|=0.592$, which holds for
$\Delta t=2^{-N}$, $N\ge1$. With $\int|a-z|\,d\mu=1.987$, the bound of
Proposition~\ref{prop:practical} at $N=12$ is $1.32\times10^{-8}$, against the
observed $6.53\times10^{-9}$. Here $\ip{H\delta_0}{\delta_0}=0$,
$\ip{H^2\delta_0}{\delta_0}=4$ and $\ip{H^3\delta_0}{\delta_0}=0$, and the leading term
of~(b), $\frac{\Delta t^2}{12}|z|$, agrees with the observed trapezoid error to three
significant digits for $N\ge6$. Subtracting the $\Delta t^2$ and $\Delta t^4$ terms of~(b)
reduces the error at $N=4$ from $4.3\times10^{-4}$ to $3.8\times10^{-7}$ and
$6.5\times10^{-10}$, and gives $1.7\times10^{-16}$ at $N=8$, from about
$5.8\times10^4$ samples.

For dissipative dynamics---absorbing boundaries, damping, open systems---$-iH$ is
replaced by the generator $B$ of a contraction semigroup, \eqref{eq:periodized} still
holds by Theorem~\ref{thm:alias}, and the rates of Theorem~\ref{thm:analytic} replace
Proposition~\ref{prop:practical}: the attainable order in $\Delta t$ is governed by the
parameter $\alpha$ of the damping, as in \S\ref{sec:num-damped}.

%% =====================================================================
\section{Discussion}\label{sec:discussion}

The dyadic resolvent identity is the telescoped rectangle rule for the Laplace
integral of the semigroup, i.e.\ the impulse-invariant discretization of the
resolvent. Its partial sums carry exactly the information of the finest sample, and
the identity is not, by itself, an efficient quadrature: when
$t\mapsto S_tx$ can be evaluated at arbitrary times, exponentially convergent rules
for the Laplace integral, or contour methods for the inverse problem, are generally
preferable~\cite{TW,WT}. What the representation does provide is an exact account of
the error of uniform sampling---the scalar $T_N/2$ and the aliases
(Corollary~\ref{cor:error})---together with the order in which successive levels
remove the aliases (Corollary~\ref{cor:levels}). This is the useful structure when the
propagator is available only on a uniform clock: autocorrelation functions
$t\mapsto\ip{S_t\psi}{\psi}$ from simulation or experiment, split-step and other
exponential integrators~\cite{HO}, and Hamiltonian simulation. In that setting the
trapezoid correction and Romberg extrapolation of Corollary~\ref{cor:romberg} turn the
first-order identity into a high-order method for smooth data, and
Theorem~\ref{thm:analytic} and Proposition~\ref{prop:analytic-sharp} say what is
achievable otherwise. Section~\ref{sec:practical} makes this explicit for self-adjoint generators. Two
practical limitations should be kept in mind. First, every
level is an infinite sum over sample times; truncating at a final time $t_{\max}$
costs a relative error of order $e^{-(\Re\lambda-\omega)t_{\max}}$, so the record must
be long compared with $1/(\Re\lambda-\omega)$, and the limit $\eps\to0^+$ of
\S\ref{sec:spectral} requires records of length $\gg1/\eps$, the usual trade-off
between spectral resolution and record length. Second, all experiments in
Section~\ref{sec:numerics} use exact samples; with noisy samples, extrapolation
trades truncation error for amplified noise, and we have not analysed this.

Several directions remain open. On Banach spaces and for $0<\alpha<1$, we do not
know whether the operator-norm rate $T^\alpha$ holds under resolvent decay
\textup{(R$_\alpha$)} alone or only up to the logarithmic factor of
Remark~\ref{rem:Ronly}. Resolvent decay slower than any power, as for semigroups that
are differentiable but not of Gevrey class, falls outside
Theorem~\ref{thm:analytic}; the bound~\eqref{eq:aliasbound} still applies on
$D(B)$ whenever $\sum_n n^{-1}\|(\lambda\pm2\pi in/T-B)^{-1}\|<\infty$. The hat
function of \S\ref{sec:num-analytic} shows that the constants in
Theorem~\ref{thm:analytic}(iii) can hide the asymptotic rate over any practical range;
a two-term bound separating low and high modes would be more informative.
Remark~\ref{rem:smooth} calls for a precise statement.
For mixed radices, Romberg extrapolation must be replaced by Neville-type schemes.
Finally, the Mellin
analogues of Remark~\ref{rem:scalar} give dyadic representations of fractional powers
$A^{-\sigma}$, spectral zeta functions and determinants, to which the alias and
Euler--Maclaurin analysis of this paper applies; these will be treated separately.

\section*{Acknowledgments}
It is a pleasure to thank Ovidiu Costin and Rodica Costin, from whose collaboration
the dyadic resolvent identity originates.

The author used Claude (Anthropic), a large language model, as an assistant in this
work: for checking proofs and computations, for writing and running the numerical
scripts, and for drafting and editing the manuscript. All mathematical statements,
proofs, code and numerical results were verified by the author, who takes full
responsibility for the content of this paper.


\begin{thebibliography}{99}

\bibitem{ABHN}
W.~Arendt, C.~J.~K.~Batty, M.~Hieber and F.~Neubrander,
\emph{Vector-valued Laplace Transforms and Cauchy Problems},
2nd ed., Monographs in Mathematics \textbf{96}, Birkh\"auser, Basel, 2011.

\bibitem{BGT}
C.~Batty, A.~Gomilko and Yu.~Tomilov,
\emph{Rational approximation of holomorphic semigroups revisited},
arXiv:2403.15894, 2024.

\bibitem{BrennerThomee}
P.~Brenner and V.~Thom\'ee,
\emph{On rational approximations of semigroups},
SIAM J. Numer. Anal. \textbf{16} (1979), 683--694.

\bibitem{CastilloThesis}
N.~Castillo,
\emph{On Rational Approximations, Resurgence and Riemann--Hilbert Problems},
Ph.D. dissertation, The Ohio State University, Columbus, OH, 2024.
OhioLINK ETD, accession no.\ osu1713373759011656.

\bibitem{CastilloV1}
N.~Castillo,
\emph{Dyadic resolvent representations of self-adjoint operators: propagator
expansions, spectral measures, and zeta functions},
arXiv:2607.21278, 2026.

\bibitem{CCC}
N.~Castillo, O.~Costin and R.~D.~Costin,
\emph{Global rational approximations of functions with factorially divergent
asymptotic series},
J. Approx. Theory \textbf{311} (2025), 106178.

\bibitem{ChenTriggiani}
S.~Chen and R.~Triggiani,
\emph{Gevrey class semigroups arising from elastic systems with gentle dissipation:
the case $0<\alpha<\frac12$},
Proc. Amer. Math. Soc. \textbf{110} (1990), 401--415.

\bibitem{CrandallPazy}
M.~G.~Crandall and A.~Pazy,
\emph{On the differentiability of weak solutions of a differential equation in
Banach space},
J. Math. Mech. \textbf{18} (1969), 1007--1016.

\bibitem{DLMF}
\emph{NIST Digital Library of Mathematical Functions}, F.~W.~J.~Olver et al., eds.,
\url{https://dlmf.nist.gov}, \S2.10.

\bibitem{Economou}
E.~N.~Economou,
\emph{Green's Functions in Quantum Physics},
3rd ed., Springer Series in Solid-State Sciences \textbf{7}, Springer, Berlin, 2006.

\bibitem{EN}
K.-J.~Engel and R.~Nagel,
\emph{One-Parameter Semigroups for Linear Evolution Equations},
Graduate Texts in Mathematics \textbf{194}, Springer, New York, 2000.

\bibitem{GomilkoTomilov}
A.~Gomilko and Yu.~Tomilov,
\emph{On convergence rates in approximation theory for operator semigroups},
J. Funct. Anal. \textbf{266} (2014), 3040--3082.

\bibitem{GKT}
A.~Gomilko, S.~Kosowicz and Yu.~Tomilov,
\emph{A general approach to approximation theory of operator semigroups},
J. Math. Pures Appl. \textbf{127} (2019), 216--267.

\bibitem{HershKato}
R.~Hersh and T.~Kato,
\emph{High-accuracy stable difference schemes for well-posed initial-value problems},
SIAM J. Numer. Anal. \textbf{16} (1979), 670--682.

\bibitem{Higham}
N.~J.~Higham,
\emph{The scaling and squaring method for the matrix exponential revisited},
SIAM J. Matrix Anal. Appl. \textbf{26} (2005), 1179--1193.

\bibitem{HO}
M.~Hochbruck and A.~Ostermann,
\emph{Exponential integrators},
Acta Numer. \textbf{19} (2010), 209--286.

\bibitem{HytonenLappas}
T.~Hyt\"onen and S.~Lappas,
\emph{Quantitative estimates for bounded holomorphic semigroups},
Semigroup Forum \textbf{108} (2024), 115--144.

\bibitem{Jackson}
L.~B.~Jackson,
\emph{A correction to impulse invariance},
IEEE Signal Process. Lett. \textbf{7} (2000), 273--275.

\bibitem{Kac}
M.~Kac,
\emph{Statistical Independence in Probability, Analysis and Number Theory},
Carus Mathematical Monographs \textbf{12}, Mathematical Association of America, 1959.

\bibitem{Kovacs}
M.~Kov\'acs,
\emph{On the convergence of rational approximations of semigroups on intermediate
spaces},
Math. Comp. \textbf{76} (2007), 273--286.

\bibitem{Lewin}
L.~Lewin,
\emph{Polylogarithms and Associated Functions},
North-Holland, New York, 1981.

\bibitem{LRT}
H.~Logemann, R.~Rebarber and S.~Townley,
\emph{Stability of infinite-dimensional sampled-data systems},
Trans. Amer. Math. Soc. \textbf{355} (2003), 3301--3328.

\bibitem{Lunardi}
A.~Lunardi,
\emph{Interpolation Theory},
3rd ed., Lecture Notes Scuola Normale Superiore \textbf{16}, Edizioni della Normale,
Pisa, 2018.

\bibitem{OS}
A.~V.~Oppenheim, R.~W.~Schafer and J.~R.~Buck,
\emph{Discrete-Time Signal Processing},
2nd ed., Prentice Hall, Upper Saddle River, NJ, 1999.

\bibitem{Pazy}
A.~Pazy,
\emph{Semigroups of Linear Operators and Applications to Partial Differential
Equations},
Applied Mathematical Sciences \textbf{44}, Springer, New York, 1983.

\bibitem{RS1}
M.~Reed and B.~Simon,
\emph{Methods of Modern Mathematical Physics, Vol.~I: Functional Analysis},
revised ed., Academic Press, New York, 1980.

\bibitem{Romberg}
W.~Romberg,
\emph{Vereinfachte numerische Integration},
Norske Vid. Selsk. Forh. (Trondheim) \textbf{28} (1955), 30--36.

\bibitem{Taylor}
S.~W.~Taylor,
\emph{Gevrey Regularity of Solutions of Evolution Equations and Boundary
Controllability},
Ph.D. thesis, University of Minnesota, 1989.

\bibitem{TE}
L.~N.~Trefethen and M.~Embree,
\emph{Spectra and Pseudospectra: The Behavior of Nonnormal Matrices and Operators},
Princeton University Press, Princeton, NJ, 2005.

\bibitem{TW}
L.~N.~Trefethen and J.~A.~C.~Weideman,
\emph{The exponentially convergent trapezoidal rule},
SIAM Rev. \textbf{56} (2014), 385--458.

\bibitem{Wakaiki24}
M.~Wakaiki,
\emph{Characterizations of the Crandall--Pazy class of $C_0$-semigroups on Hilbert
spaces and their application to decay estimates},
J. Funct. Anal. \textbf{289} (2025), no.~3, 110902.

\bibitem{Wakaiki25}
M.~Wakaiki,
\emph{Relation between semigroup growth and resolvent decay for immediately
differentiable semigroups},
J. Evol. Equ. \textbf{26} (2026), 79.

\bibitem{WT}
J.~A.~C.~Weideman and L.~N.~Trefethen,
\emph{Parabolic and hyperbolic contours for computing the Bromwich integral},
Math. Comp. \textbf{76} (2007), 1341--1356.

\end{thebibliography}
\end{document}